\documentclass[final,3p]{elsarticle}
\usepackage{eurosym}
\usepackage{amsfonts}
\usepackage{amsmath}
\usepackage{amssymb}
\usepackage{graphicx}
\usepackage[tableposition=top]{caption}
\usepackage{subcaption}
\usepackage{float}
\usepackage{setspace}
\usepackage{placeins}
\usepackage{mathrsfs}
\usepackage{color}
\usepackage{booktabs}
\usepackage{hyperref}

\newsavebox \foobox
\newlength{\foodim}

\graphicspath{{C:/Figures/}}

\newtheorem{definition}{Definition}
\newtheorem{example}{Example}

\newtheorem{remark}{Remark}

\journal{ }

\begin{document}

\begin{frontmatter}
\title{A Newton interpolation based predictor-corrector numerical method for fractional differential equations with an activator-inhibitor case study}
\author{Redouane Douaifia$^{a}$, Samir Bendoukha$^{b}$و Salem Abdelmalek$^{a,c,*}$}
\address{(a) Laboratory of Mathematics, Informatics and Systems (LAMIS), Larbi Tebessi University - Tebessa, Algeria\\
	(b) Electrical Engineering Department, College of Engineering at Yanbu, Taibah University, Saudi Arabia\\ 
	(c) Department of Mathematics and Computer Science, Larbi Tebessi University - Tebessa, Algeria\\
	(*) Corresponding author, email: salem.abdelmalek@univ-tebessa.dz}
\begin{abstract}
This paper presents a new predictor-corrector numerical scheme suitable for fractional differential equations. An improved explicit Atangana-Seda formula is obtained by considering the neglected terms and used as the predictor stage of the proposed method. Numerical formulas are presented that approximate the classical first derivative as well as the Caputo, Caputo-Fabrizio and Atangana-Baleanu fractional derivatives. Simulation results are used to assess the approximation error of the new method for various differential equations. In addition, a case study is considered where the proposed scheme is used to obtained numerical solutions of the Gierer-Meinhardt activator-inhibitor model with the aim of assessing the system's dynamics.
\end{abstract}
\end{frontmatter}

\section{Introduction}

Over the last century, ordinary and partial differential equations have been
shown to produce accurate models of real life phenomena spanning a range of
different scientific and engineering disciplines. Based on these models,
researchers are able to infer the characteristics of these phenomena and
devise effective control strategies. Such characteristics include the
existence and boundedness of solutions, blow-up time, asymptotic behavior,
and more. Since these models can be quite complicated and analytical
solutions are not always attainable, numerical analysis became a useful tool
that helps obtain approximate solutions and give indications on the behavior
of these models. The simplest numerical methods reported in the literature
and suitable for linear systems are based on linear interpolation, which has
been around for over 2000 years. For the nonlinear case, well established
interpolation techniques include Newton's method, Lagrange interpolation
polynomials, Gaussian elimination, and Euler's method \cite%
{Werner1984,Fred1970,Yang2015,Dimitrov1994}.

In recent years, an apparent shift has been observed from classic models
involving integer-order derivatives to fractional ones. This shift may be
attributed to the many benefits associated with fractional derivatives
including their infinite memory and wider dynamical range. Numerical methods
had to evolve in order for researchers to investigate these fractional
models. Several numerical schemes have been proposed for solving fractional
ordinary differential equations, especially nonlinear ones including \cite%
{Diethelm1998,Ford2001,Odibat2008,Moghaddam2016,Asl2017,Patricio2019}%
. To the best of the authors' knowledge, the most widely accepted scheme is
the Adams-Bashforth method developed with a Lagrange interpolation
polynomial basis \cite{Zhang2018,Jain2018}. In recent years, studies have
shown that on average, Newton's method is superior to Lagrange polynomials
taking into consideration a wide range of polynomial functions \cite%
{Srivastava2012a,Srivastava2012b}. A numerical method suitable for both
integer and fractional ordinary differential systems was proposed by
Atangana and Seda by replacing the Lagrange polynomial interpolation of the
Adams-Bashforth scheme with Newton quadratic interpolation in \cite%
{Atangana2020,Atangana2020Corrigendum}. The authors derived iterative
numerical formulas for the standard and fractal versions of the Caputo,
Caputo-Fabrizio, and Atangana-Baleanu fractional derivatives. This method
was applied to chaotic systems and showed promising results \cite%
{Alkahtani2019,Atangana2020a,Atangana2020c}. The method was also extended to
partial differential equations with integer and non-integer orders \cite%
{Atangana2020b}.

Over the last few decades a class of numerical methods called
predictor-corrector emerged and became the center of attention for many
researchers \cite{Gragg1964,Marciniak2017,Butcher2016}. It is well known
that numerical methods are generally divided into implicit and explicit
types and that the implicit type is more stable and efficient but difficult
to solve due to the fact that the unknown appears on both sides of the
formula. Predictor-corrector methods work in two steps. An initial explicit
approximation (predictor) of the solution is obtained and substituted into
right side of the implicit formula (corrector). A predictor-corrector
Adams-Bashforth method was introduced in \cite{Diethelm2002}. In this
method, the explicit one-step Adams--Bashforth rule and the implicit
one-step Adams-Moulton method are used as predictor and corrector,
respectively. Other more recent works include \cite%
{Nguyen2017,Douaifia2019,Kumar2019,Heris2019}. In this paper, we propose a
new predictor-corrector method where an improved version of the
Atangana-Seda method of \cite{Atangana2020,Atangana2020Corrigendum} is used
as the predictor. We derive iterative formulas for the classical as well as
the Caputo, Caputo-Fabrizio and Atangana-Baleanu fractional derivative
scenarios. Numerical examples are presented to evaluate the effectiveness of
the proposed methods.

\section{Important Definitions}

Before we delve into the main concern of the paper, let us describe the
fractional integrals and derivatives that will be used in our work. For more
on these definitions, the reader may wish to refer to \cite%
{Atangana2020,Podlubny,Kilbas2006,CaputoFabrizio2015,AtanganaBaleanu2016}.

\begin{definition}
\label{Def1}The $\alpha $--order Riemann--Liouville fractional integral of a
function $x(t)$ is defined as 
\begin{equation}
_{0}I_{t}^{\alpha }x\left( t\right) =\frac{1}{\Gamma \left( \alpha \right) }%
\int_{0}^{t}\left( t-s \right) ^{\alpha -1}x(s)ds,  \label{1.1}
\end{equation}%
where $\alpha >0$ and $\Gamma \left( \alpha \right) $\ is the Gamma function
defined as%
\begin{equation}
\Gamma \left( \alpha \right) =\int_{0}^{+\infty }e^{-t}t^{\alpha -1}dt,
\label{1.2}
\end{equation}%
for $Re(\alpha )>0$.
\end{definition}

\begin{definition}
\label{Def2}The $\alpha $--order Caputo fractional derivative of a function $%
x(t)$ is defined as%
\begin{equation}
_{0}^{C}D_{t}^{\alpha }x\left( t\right) =\left \{ 
\begin{array}{l}
_{0}I_{t}^{n-\alpha }\left \{ \frac{d^{n}}{dt^{n}}x\left( t\right) \right \} 
\text{, if }n-1<\alpha <n\in \mathbb{N,\medskip } \\ 
\frac{d^{n}}{dt^{n}}x\left( t\right) \text{, if }\alpha =n\in \mathbb{N}.%
\end{array}%
\right.  \label{1.3}
\end{equation}
\end{definition}

\begin{definition}
\label{Def3}The Caputo-Fabrizio fractional integral of a function $x(t)$ is
defined as 
\begin{equation}
_{0}^{CF}I_{t}^{\alpha }x\left( t\right) =\frac{1-\alpha }{M(\alpha )}x(t)+%
\frac{\alpha }{M(\alpha )}\int_{0}^{t}x(s)ds,  \label{1.4}
\end{equation}%
where $\alpha \in (0,1)$, and $M(\alpha )$ is a normalization function
satisfying $M(0)=M(1)=1$.
\end{definition}

\begin{definition}
\label{Def4}Let $x\in H^{1}\left( [0,T]\right) $, $T>0$, and $\alpha \in
(0,1)$. The Caputo-Fabrizio fractional derivative of a function $x(t)$ is
defined as%
\begin{equation}
_{0}^{CF}D_{t}^{\alpha }x\left( t\right) =\frac{M(\alpha )}{1-\alpha }%
\int_{0}^{t}\frac{d}{ds}x\left( s\right) \exp \left( -\frac{\alpha (t-s)}{%
1-\alpha }\right) ds.  \label{1.5}
\end{equation}
\end{definition}

\begin{definition}
\label{Def5}The Atangana-Baleanu fractional integral of a function $x(t)$ is
defined as%
\begin{equation}
_{0}^{ABC}I_{t}^{\alpha }x\left( t\right) =\frac{1-\alpha }{AB(\alpha )}x(t)+%
\frac{\alpha }{AB(\alpha )\Gamma (\alpha )}\int_{0}^{t}x(s)\left( t-s \right) ^{\alpha -1}ds,  \label{1.6}
\end{equation}%
where $\alpha \in (0,1)$, and 
\begin{equation}
AB(\alpha )=1-\alpha +\frac{\alpha }{\Gamma (\alpha )}.  \label{1.7}
\end{equation}
\end{definition}

\begin{definition}
\label{Def6}Let $x\in H^{1}\left( [0,T]\right) $, $T>0$, and $\alpha \in
(0,1)$. The Atangana-Baleanu fractional derivative in the\ Caputo sense of a
function $x(t)$ is defined as%
\begin{equation}
_{0}^{ABC}D_{t}^{\alpha }x\left( t\right) =\frac{AB(\alpha )}{1-\alpha }%
\int_{0}^{t}\frac{d}{ds}x\left( s\right) E_{\alpha }\left( -\frac{\alpha
(t-s)^{\alpha }}{1-\alpha }\right) ds,  \label{1.8}
\end{equation}%
where $E_{\alpha }(z)$ is the Mittag-Leffler kernel function of order $%
\alpha $ defined as%
\begin{equation}
E_{\alpha }(z)=\sum_{k=0}^{\infty }\frac{z^{k}}{\Gamma (\alpha k+1)},
\label{1.9}
\end{equation}%
for $Re(\alpha )>0$ and$\ z\in \mathbb{C}$.
\end{definition}

\section{The Proposed Predictor-Corrector Method}

\subsection{Classical Derivative}

We start with the simple classical initial-value problem given by%
\begin{equation}
\left \{ 
\begin{array}{l}
\frac{dy(t)}{dt}=f(t,y(t)),\mathbb{\medskip } \\ 
y(0)=y_{0},%
\end{array}%
\right.  \label{2.1}
\end{equation}%
where $f$ is a smooth nonlinear function guaranteeing a unique solution $%
y(t) $. In order to develop a numerical formula approximating the solution
of (\ref{2.1}), we convert the differential equation into the integral%
\begin{equation}
y(t)-y(0)=\int_{0}^{t}f(s,y(s))ds.  \label{2.2}
\end{equation}%
In an iterative approximation, we may choose two distinct points in time $%
t_{m}=m\Delta t$\ and $t_{m+1}=(m+1)\Delta t$. Substituting these points
into (\ref{2.2}) yields%
\begin{equation*}
y\left( t_{m}\right) -y(0)=\int_{0}^{t_{m}}f(s,y(s))ds,
\end{equation*}%
and%
\begin{equation*}
y\left( t_{m+1}\right) -y(0)=\int_{0}^{t_{m+1}}f(s,y(s))ds,
\end{equation*}%
respectively. Taking the difference yields%
\begin{equation}
y\left( t_{m+1}\right) -y\left( t_{m}\right)
=\int_{t_{m}}^{t_{m+1}}f(s,y(s))ds.  \label{2.3}
\end{equation}%
Hence, the function $f(s,y(s))$\ may be approximated over the interval $%
[t_{m},t_{m+1}]$\ by means of Newton's second order interpolation polynomial
given by%
\begin{eqnarray}
\mathcal{N}_{m}(s) &=&f\left( t_{m+1},y\left( t_{m+1}\right) \right) +\frac{%
f\left( t_{m+1},y\left( t_{m+1}\right) \right) -f\left( t_{m},y\left(
t_{m}\right) \right) }{\Delta t}\left( s-t_{m+1}\right)  \notag \\
&&\  \  \  \ +\frac{f\left( t_{m+1},y\left( t_{m+1}\right) \right) -2f\left(
t_{m},y\left( t_{m}\right) \right) +f\left( t_{m-1},y\left( t_{m-1}\right)
\right) }{2(\Delta t)^{2}}  \notag \\
&&\  \  \  \  \times \left( s-t_{m}\right) \left( s-t_{m+1}\right) .  \label{2.4}
\end{eqnarray}%
Substitution into (\ref{2.3}) leads to the difference formula%
\begin{eqnarray}
y_{m+1}-y_{m} &=&f\left( t_{m+1},y_{m+1}\right) \Delta t+\left( \frac{%
f\left( t_{m+1},y_{m+1}\right) -f\left( t_{m},y_{m}\right) }{\Delta t}%
\right) \int_{t_{m}}^{t_{m+1}}\left( s-t_{m+1}\right) ds  \notag \\
&&\  \  \  \ +\left( \frac{f\left( t_{m+1},y_{m+1}\right) -2f\left(
t_{m},y_{m}\right) +f\left( t_{m-1},y_{m-1}\right) }{2(\Delta t)^{2}}\right)
\notag \\
&&\  \  \  \  \  \  \  \  \int_{t_{m}}^{t_{m+1}}\left( s-t_{m}\right) \left(
s-t_{m+1}\right) ds.  \label{2.5}
\end{eqnarray}%
Given that%
\begin{equation}
\int_{t_{m}}^{t_{m+1}}\left( s-t_{m+1}\right) ds=-\frac{(\Delta t)^{2}}{2},
\label{2.5.1}
\end{equation}%
and%
\begin{equation}
\int_{t_{m}}^{t_{m+1}}\left( s-t_{m}\right) \left( s-t_{m+1}\right) ds=-%
\frac{(\Delta t)^{3}}{6},  \label{2.5.2}
\end{equation}%
formula (\ref{2.5}) reduces to the implicit form%
\begin{eqnarray}
y_{m+1}-y_{m} &=&f\left( t_{m+1},y_{m+1}\right) \Delta t-\left[ f\left(
t_{m+1},y_{m+1}\right) -f\left( t_{m},y_{m}\right) \right] \frac{\Delta t}{2}
\notag \\
&&-\left[ f\left( t_{m+1},y_{m+1}\right) -2f\left( t_{m},y_{m}\right)
+f\left( t_{m-1},y_{m-1}\right) \right] \frac{\Delta t}{12}.  \label{2.6}
\end{eqnarray}%
The term $y_{m+1}$ appears on both sides of the formula. The
predictor-corrector scheme works by first producing an approximation of $%
y_{m+1}$ denoted by $y_{m+1}^{P}$, and then using (\ref{2.6}) to correct the
approximation. The correction formula is, thus, given by%
\begin{equation}
y_{m+1}=y_{m}+\frac{5}{12}f\left( t_{m+1},y_{m+1}^{P}\right) \Delta t+\frac{2%
}{3}f\left( t_{m},y_{m}\right) \Delta t-f\left( t_{m-1},y_{m-1}\right) \frac{%
\Delta t}{12},  \label{2.7}
\end{equation}%
where the predictor $y_{m+1}^{P}$ is obtained by means of the Atangana-Seda
scheme (cf. \cite{Atangana2020}), i.e.%
\begin{equation}
y_{m+1}^{P}=y_{m}+\frac{5}{12}f\left( t_{m-2},y_{m-2}\right) \Delta t-\frac{4%
}{3}f\left( t_{m-1},y_{m-1}\right) \Delta t+\frac{23}{12}f\left(
t_{m},y_{m}\right) \Delta t.  \label{2.8}
\end{equation}

\subsection{Caputo Fractional Derivative\label{Section_Caputo}}

Let us now move to the fractional derivative case. Various derivatives have
been proposed throughout the years. However, the most commonly used is the
Caputo one. We consider the initial-value problem%
\begin{equation}
\left \{ 
\begin{array}{l}
_{0}^{C}D_{t}^{\alpha }y(t)=f(t,y(t)),\mathbb{\medskip } \\ 
y(0)=y_{0},%
\end{array}%
\right.  \label{3.1}
\end{equation}%
with $\alpha \in (0,1]$, and $f$ being a smooth nonlinear function such that
(\ref{3.1}) admits a unique solution $y(t)$. Following the same procedure of
the standard case, we start with the integral%
\begin{equation}
y(t)-y(0)=\frac{1}{\Gamma (\alpha )}\int_{0}^{t}f(s,y(s))(t-s)^{\alpha -1}ds.
\label{3.2}
\end{equation}%
At the single point $t_{m+1}=(m+1)\Delta t$, we have the following 
\begin{eqnarray}
y\left( t_{m+1}\right) &=&y(0)+\frac{1}{\Gamma (\alpha )}%
\int_{0}^{t_{m+1}}f(s,y(s))\left( t_{m+1}-s\right) ^{\alpha -1}ds  \notag \\
&=&y(0)+\frac{1}{\Gamma (\alpha )}\sum_{i=0}^{m}%
\int_{t_{i}}^{t_{i+1}}f(s,y(s))\left( t_{m+1}-s\right) ^{\alpha -1}ds,
\label{3.3}
\end{eqnarray}%
with $t_{0}=0$. Function $f(s,y(s))$ can be approximated over the
sub-interval $[t_{i},t_{i+1}]$ as a polynomial by means of%
\begin{equation}
\mathcal{N}_{i}(s)=\left \{ 
\begin{array}{lll}
\widetilde{\mathcal{N}}_{i}(s) & \text{if} & i=0,\mathbb{\medskip } \\ 
\widehat{\mathcal{N}}_{i}(s) & \text{if} & i\in \left \{ 1,\dots ,m\right \}
,%
\end{array}%
\right.  \label{3.4}
\end{equation}%
where%
\begin{equation}
\widetilde{\mathcal{N}}_{i}(s)=f(t_{i},y(t_{i}))+\left( \frac{%
f(t_{i+1},y(t_{i+1}))-f(t_{i},y(t_{i}))}{\Delta t}\right) (s-t_{i}),
\label{3.5}
\end{equation}%
and%
\begin{eqnarray}
\widehat{\mathcal{N}}_{i}(s) &=&f\left( t_{i+1},y\left( t_{i+1}\right)
\right) +\frac{f\left( t_{i+1},y\left( t_{i+1}\right) \right) -f\left(
t_{i},y\left( t_{i}\right) \right) }{\Delta t}\left( s-t_{i+1}\right)  \notag
\\
&&+\frac{f\left( t_{i+1},y\left( t_{i+1}\right) \right) -2f\left(
t_{i},y\left( t_{i}\right) \right) +f\left( t_{i-1},y\left( t_{i-1}\right)
\right) }{2(\Delta t)^{2}}\times \left( s-t_{i}\right) \left(
s-t_{i+1}\right) .  \label{3.6}
\end{eqnarray}%
Using the Newton polynomial (\ref{3.4}), formula (\ref{3.3}) becomes%
\begin{eqnarray}
y(t_{m+1}) &=&y(0)+\frac{1}{\Gamma (\alpha )}\int_{0}^{t_{1}}\left[
f(t_{0},y(t_{0}))+\left( \frac{f(t_{1},y(t_{1}))-f(t_{0},y(t_{0}))}{\Delta t}%
\right) s\right] \left( t_{m+1}-s\right) ^{\alpha -1}ds  \notag \\
&&+\frac{1}{\Gamma (\alpha )}\sum_{i=1}^{m}\int_{t_{i}}^{t_{i+1}}\left \{ 
\begin{array}{l}
f\left( t_{i+1},y\left( t_{i+1}\right) \right) \\ 
+\frac{f\left( t_{i+1},y\left( t_{i+1}\right) \right) -f\left( t_{i},y\left(
t_{i}\right) \right) }{\Delta t}\left( s-t_{i+1}\right) \\ 
+\frac{f\left( t_{i+1},y\left( t_{i+1}\right) \right) -2f\left(
t_{i},y\left( t_{i}\right) \right) +f\left( t_{i-1},y\left( t_{i-1}\right)
\right) }{2(\Delta t)^{2}} \\ 
\times \left( s-t_{i}\right) \left( s-t_{i+1}\right)%
\end{array}%
\right \} \left( t_{m+1}-s\right) ^{\alpha -1}ds.  \label{3.7}
\end{eqnarray}%
Simplifying and rearranging the terms leads to%
\begin{eqnarray}
y_{m+1} &=&y_{0}+\frac{1}{\Gamma (\alpha )}f(t_{0},y_{0})\int_{0}^{t_{1}}%
\left( t_{m+1}-s\right) ^{\alpha -1}ds+\frac{1}{\Gamma (\alpha )}\left( 
\frac{f(t_{1},y_{1})-f(t_{0},y_{0})}{\Delta t}\right)
\int_{0}^{t_{1}}s\left( t_{m+1}-s\right) ^{\alpha -1}ds  \notag \\
&&+\frac{1}{\Gamma (\alpha )}\sum_{i=1}^{m}f\left( t_{i+1},y_{i+1}\right)
\int_{t_{i}}^{t_{i+1}}\left( t_{m+1}-s\right) ^{\alpha -1}ds  \notag \\
&&+\frac{1}{\Gamma (\alpha )}\sum_{i=1}^{m}\frac{f\left(
t_{i+1},y_{i+1}\right) -f\left( t_{i},y_{i}\right) }{\Delta t}%
\int_{t_{i}}^{t_{i+1}}\left( s-t_{i+1}\right) \left( t_{m+1}-s\right)
^{\alpha -1}ds  \notag \\
&&+\frac{1}{\Gamma (\alpha )}\sum_{i=1}^{m}\frac{f\left(
t_{i+1},y_{i+1}\right) -2f\left( t_{i},y_{i}\right) +f\left(
t_{i-1},y_{i-1}\right) }{2(\Delta t)^{2}}  \notag \\
&&\times \int_{t_{i}}^{t_{i+1}}\left( s-t_{i}\right) \left( s-t_{i+1}\right)
\left( t_{m+1}-s\right) ^{\alpha -1}ds  \label{3.8}
\end{eqnarray}%
The four different integrals in (\ref{3.8}) can be calculated as%
\begin{equation}
\int_{0}^{t_{1}}s\left( t_{m+1}-s\right) ^{\alpha -1}ds=\frac{(\Delta
t)^{\alpha +1}}{\alpha (\alpha +1)}\left[ (m+1)^{\alpha +1}-m^{\alpha
+1}-(\alpha +1)m^{\alpha }\right] ,  \label{3.9}
\end{equation}%
\begin{equation}
\int_{t_{i}}^{t_{i+1}}\left( t_{m+1}-s\right) ^{\alpha -1}ds=\frac{(\Delta
t)^{\alpha }}{\alpha }\left[ (m-i+1)^{\alpha }-(m-i)^{\alpha }\right] ,
\label{3.10}
\end{equation}%
\begin{equation}
\int_{t_{i}}^{t_{i+1}}\left( s-t_{i+1}\right) \left( t_{m+1}-s\right)
^{\alpha -1}ds=\frac{(\Delta t)^{\alpha +1}}{\alpha (\alpha +1)}\left[
(m-i-\alpha )(m-i+1)^{\alpha }-(m-i)^{\alpha +1}\right] ,  \label{3.11}
\end{equation}%
and%
\begin{eqnarray}
\int_{t_{i}}^{t_{i+1}}\left( s-t_{i}\right) \left( s-t_{i+1}\right) \left(
t_{m+1}-s\right) ^{\alpha -1}ds &=&\frac{(\Delta t)^{\alpha +2}}{\alpha
(\alpha +1)(\alpha +2)}  \notag \\
&&\times \left[ 
\begin{array}{l}
(m-i+1)^{\alpha }\left[ 
\begin{array}{c}
2(m-i)^{2}-\alpha (m-i+1) \\ 
+2(m-i)%
\end{array}%
\right] \\ 
-(m-i)^{\alpha }\left[ 
\begin{array}{c}
2(m-i)^{2}+\alpha (m-i) \\ 
+2(m-i)%
\end{array}%
\right]%
\end{array}%
\right] ,  \label{3.12}
\end{eqnarray}%
respectively. By substituting these calculations into (\ref{3.8}), we
obtain%
\begin{eqnarray}
y_{m+1} &=&y_{0}+\frac{(\Delta t)^{\alpha }}{\Gamma (\alpha +1)}%
f(t_{0},y_{0})\left[ (m+1)^{\alpha }-m^{\alpha }\right]  \notag \\
&&+\frac{(\Delta t)^{\alpha }}{\Gamma (\alpha +2)}\left(
f(t_{1},y_{1})-f(t_{0},y_{0})\right) \left[ (m+1)^{\alpha +1}-m^{\alpha
+1}-(\alpha +1)m^{\alpha }\right]  \notag \\
&&+\frac{(\Delta t)^{\alpha }}{\Gamma (\alpha +1)}\sum_{i=1}^{m}f\left(
t_{i+1},y_{i+1}\right) \left[ (m-i+1)^{\alpha }-(m-i)^{\alpha }\right] 
\notag \\
&&+\frac{(\Delta t)^{\alpha }}{\Gamma (\alpha +2)}\sum_{i=1}^{m}\left(
f\left( t_{i+1},y_{i+1}\right) -f\left( t_{i},y_{i}\right) \right) \left[
(m-i-\alpha )(m-i+1)^{\alpha }-(m-i)^{\alpha +1}\right]  \notag \\
&&+\frac{(\Delta t)^{\alpha }}{2\Gamma (\alpha +3)}\sum_{i=1}^{m}\left(
f\left( t_{i+1},y_{i+1}\right) -2f\left( t_{i},y_{i}\right) +f\left(
t_{i-1},y_{i-1}\right) \right)  \notag \\
&&\times \left[ 
\begin{array}{l}
(m-i+1)^{\alpha }\left[ 
\begin{array}{c}
2(m-i)^{2}-\alpha (m-i+1) \\ 
+2(m-i)%
\end{array}%
\right] \\ 
-(m-i)^{\alpha }\left[ 
\begin{array}{c}
2(m-i)^{2}+\alpha (m-i) \\ 
+2(m-i)%
\end{array}%
\right]%
\end{array}%
\right] .  \label{3.13}
\end{eqnarray}%
In order to simplify the formulas to come, let us define the expresion%
\begin{eqnarray}
\Upsilon _{p} &=&\frac{(\Delta t)^{\alpha }}{\Gamma (\alpha +1)}%
\sum_{i=1}^{p}f\left( t_{i+1},y_{i+1}\right) \left[ (m-i+1)^{\alpha
}-(m-i)^{\alpha }\right]  \notag \\
&&+\frac{(\Delta t)^{\alpha }}{\Gamma (\alpha +2)}\sum_{i=1}^{p}\left(
f\left( t_{i+1},y_{i+1}\right) -f\left( t_{i},y_{i}\right) \right) \left[
(m-i-\alpha )(m-i+1)^{\alpha }-(m-i)^{\alpha +1}\right]  \notag \\
&&+\frac{(\Delta t)^{\alpha }}{2\Gamma (\alpha +3)}\sum_{i=1}^{p}\left(
f\left( t_{i+1},y_{i+1}\right) -2f\left( t_{i},y_{i}\right) +f\left(
t_{i-1},y_{i-1}\right) \right)  \notag \\
&&\times \left[ 
\begin{array}{l}
(m-i+1)^{\alpha }\left[ 
\begin{array}{c}
2(m-i)^{2}-\alpha (m-i+1) \\ 
+2(m-i)%
\end{array}%
\right] \\ 
-(m-i)^{\alpha }\left[ 
\begin{array}{c}
2(m-i)^{2}+\alpha (m-i) \\ 
+2(m-i)%
\end{array}%
\right]%
\end{array}%
\right] ,  \label{3.14}
\end{eqnarray}%
with the convention%
\begin{equation}
\Upsilon _{0}=0.  \label{3.15}
\end{equation}%
Using this notation, (\ref{3.13}) can be rewritten in the form%
\begin{eqnarray}
y_{m+1} &=&y_{0}+\Upsilon _{m-1}+\frac{(\Delta t)^{\alpha }}{\Gamma (\alpha
+1)}f(t_{0},y_{0})\left[ (m+1)^{\alpha }-m^{\alpha }\right]  \notag \\
&&+\frac{(\Delta t)^{\alpha }}{\Gamma (\alpha +2)}\left(
f(t_{1},y_{1})-f(t_{0},y_{0})\right) \left[ (m+1)^{\alpha +1}-m^{\alpha
+1}-(\alpha +1)m^{\alpha }\right]  \notag \\
&&+\frac{(\Delta t)^{\alpha }}{\Gamma (\alpha +1)}f\left(
t_{m+1},y_{m+1}\right) +\frac{\alpha (\Delta t)^{\alpha }}{\Gamma (\alpha +2)%
}\left( f\left( t_{m},y_{m}\right) -f\left( t_{m+1},y_{m+1}\right) \right) 
\notag \\
&&-\frac{\alpha (\Delta t)^{\alpha }}{2\Gamma (\alpha +3)}\left( f\left(
t_{m+1},y_{m+1}\right) -2f\left( t_{m},y_{m}\right) +f\left(
t_{m-1},y_{m-1}\right) \right) .  \label{3.16}
\end{eqnarray}

Formula (\ref{3.16}) will serve as our implicit part, i.e. the corrector.
The terms $y_{m+1}$ on the right hand side will be replaced by the predictor 
$y_{m+1}^{P}$, which will be an improved version of the Atangana-Seda scheme
derived for the Caputo fractional derivative in \cite{Atangana2020}. To
obtain our predictor formula, let us go back to (\ref{3.2}) and use the
predictor notation $y^{P}(t)$, which yields%
\begin{equation*}
y^{P}(t)-y(0)=\frac{1}{\Gamma (\alpha )}\int_{0}^{t}f(s,y(s))(t-s)^{\alpha
-1}ds,
\end{equation*}%
and, consequently, at $t_{m+1}=(m+1)\Delta t$, we have 
\begin{equation}
y^{P}\left( t_{m+1}\right) =y(0)+\frac{1}{\Gamma (\alpha )}%
\sum_{i=0}^{m}\int_{t_{i}}^{t_{i+1}}f(s,y(s))\left( t_{m+1}-s\right)
^{\alpha -1}ds.  \label{3.17}
\end{equation}%
The function $f(s,y(s))$ can be approximated over each sub-interval $%
[t_{i},t_{i+1}]$ using a delayed version of the Newton's
polynomial seen earlier in (\ref{3.4}) and given by%
\begin{equation}
\mathcal{N}_{i}^{P}(s)=\left \{ 
\begin{array}{lll}
\widetilde{\mathcal{N}}_{i}^{P}(s) & \text{if} & i\in \left \{ 0,1\right \}
,\medskip \\ 
\widehat{\mathcal{N}}_{i}^{P}(s) & \text{if} & i\in \left \{ 2,\dots
,m\right \} ,%
\end{array}%
\right.  \label{3.18}
\end{equation}%
where%
\begin{equation}
\widetilde{\mathcal{N}}_{i}^{P}(s)=f(t_{i},y(t_{i}))+\left( \frac{%
f(t_{i+1},y(t_{i+1}))-f(t_{i},y(t_{i}))}{\Delta t}\right) (s-t_{i}),
\label{3.19}
\end{equation}%
and%
\begin{eqnarray}
\widehat{\mathcal{N}}_{i}^{P}(s) &=&f\left( t_{i-2},y\left( t_{i-2}\right)
\right) +\frac{f\left( t_{i-1},y\left( t_{i-1}\right) \right) -f\left(
t_{i-2},y\left( t_{i-2}\right) \right) }{\Delta t}\left( s-t_{i-2}\right) 
\notag \\
&&+\frac{f\left( t_{i},y\left( t_{i}\right) \right) -2f\left(
t_{i-1},y\left( t_{i-1}\right) \right) +f\left( t_{i-2},y\left(
t_{i-2}\right) \right) }{2(\Delta t)^{2}}\times \left( s-t_{i-2}\right)
\left( s-t_{i-1}\right) .  \label{3.20}
\end{eqnarray}%
Substituting the interpolated approximation of $f(s,y(s))$ into (\ref{3.17})
yields the predictor%
\begin{eqnarray}
y_{m+1}^{P} &=&y_{0}+\frac{1}{\Gamma (\alpha )}\sum_{i=0}^{1}f\left(
t_{i},y_{i}\right) \int_{t_{i}}^{t_{i+1}}\left( t_{m+1}-s\right) ^{\alpha
-1}ds  \notag \\
&&+\frac{1}{\Gamma (\alpha )}\sum_{i=0}^{1}\frac{f\left(
t_{i+1},y_{i+1}\right) -f\left( t_{i},y_{i}\right) }{\Delta t}%
\int_{t_{i}}^{t_{i+1}}\left( s-t_{i}\right) \left( t_{m+1}-s\right) ^{\alpha
-1}ds  \notag \\
&&+\frac{1}{\Gamma (\alpha )}\sum_{i=2}^{m}f\left( t_{i-2},y_{i-2}\right)
\int_{t_{i}}^{t_{i+1}}\left( t_{m+1}-s\right) ^{\alpha -1}ds  \notag \\
&&+\frac{1}{\Gamma (\alpha )}\sum_{i=2}^{m}\frac{f\left(
t_{i-1},y_{i-1}\right) -f\left( t_{i-2},y_{i-2}\right) }{\Delta t}%
\int_{t_{i}}^{t_{i+1}}\left( s-t_{i-2}\right) \left( t_{m+1}-s\right)
^{\alpha -1}ds  \notag \\
&&+\frac{1}{\Gamma (\alpha )}\sum_{i=2}^{m}\frac{f\left( t_{i},y_{i}\right)
-2f\left( t_{i-1},y_{i-1}\right) +f\left( t_{i-2},y_{i-2}\right) }{2(\Delta
t)^{2}}  \notag \\
&&\times \int_{t_{i}}^{t_{i+1}}\left( s-t_{i-2}\right) \left(
s-t_{i-1}\right) \left( t_{m+1}-s\right) ^{\alpha -1}ds.  \label{3.21}
\end{eqnarray}%
We can calculate the integrals as%
\begin{equation}
\int_{t_{i}}^{t_{i+1}}\left( s-t_{i}\right) \left( t_{m+1}-s\right) ^{\alpha
-1}ds=\frac{(\Delta t)^{\alpha +1}}{\alpha (\alpha +1)}\left[
(m-i+1)^{\alpha +1}-(m-i)^{\alpha +1}-(\alpha +1)(m-i)^{\alpha }\right] ,
\label{3.22}
\end{equation}%
\begin{equation}
\int_{t_{i}}^{t_{i+1}}\left( s-t_{i-2}\right) \left( t_{m+1}-s\right)
^{\alpha -1}ds=\frac{(\Delta t)^{\alpha +1}}{\alpha (\alpha +1)}\left[ 
\begin{array}{c}
(m-i+1)^{\alpha }(m-i+3+2\alpha ) \\ 
-(m-i)^{\alpha }(m-i+3+3\alpha )%
\end{array}%
\right] ,  \label{3.23}
\end{equation}%
and%
\begin{eqnarray}
\int_{t_{i}}^{t_{i+1}}\left( s-t_{i-2}\right) \left( s-t_{i-1}\right) \left(
t_{m+1}-s\right) ^{\alpha -1}ds &=&\left[ 
\begin{array}{c}
(m-i+1)^{\alpha }\left[ 
\begin{array}{c}
2(m-i)^{2}+(3\alpha +10)(m-i) \\ 
+2\alpha ^{2}+9\alpha +12%
\end{array}%
\right] \\ 
-(m-i)^{\alpha }\left[ 
\begin{array}{c}
2(m-i)^{2}+(5\alpha +10)(m-i) \\ 
+6\alpha ^{2}+18\alpha +12%
\end{array}%
\right]%
\end{array}%
\right]  \notag \\
&&\times \frac{(\Delta t)^{\alpha +2}}{\alpha (\alpha +1)(\alpha +2)}.
\label{3.24}
\end{eqnarray}%
Substituting these calculations into (\ref{3.21}) produces the improved
Atangana-Seda scheme predictor%
\begin{eqnarray}
y_{m+1}^{P} &=&y_{0}+\frac{(\Delta t)^{\alpha }}{\Gamma (\alpha +1)}%
\sum_{i=0}^{1}f\left( t_{i},y_{i}\right) \left[ (m-i+1)^{\alpha
}-(m-i)^{\alpha }\right]  \notag \\
&&+\frac{(\Delta t)^{\alpha }}{\Gamma (\alpha +2)}\sum_{i=0}^{1}\left(
f\left( t_{i+1},y_{i+1}\right) -f\left( t_{i},y_{i}\right) \right)  \notag \\
&&\times \left[ (m-i+1)^{\alpha +1}-(m-i)^{\alpha +1}-(\alpha
+1)(m-i)^{\alpha }\right]  \notag \\
&&+\frac{(\Delta t)^{\alpha }}{\Gamma (\alpha +1)}\sum_{i=2}^{m}f\left(
t_{i-2},y_{i-2}\right) \left[ (m-i+1)^{\alpha }-(m-i)^{\alpha }\right] 
\notag \\
&&+\frac{(\Delta t)^{\alpha }}{\Gamma (\alpha +2)}\sum_{i=2}^{m}\left(
f\left( t_{i-1},y_{i-1}\right) -f\left( t_{i-2},y_{i-2}\right) \right) 
\notag \\
&&\times \left[ 
\begin{array}{c}
(m-i+1)^{\alpha }(m-i+3+2\alpha ) \\ 
-(m-i)^{\alpha }(m-i+3+3\alpha )%
\end{array}%
\right]  \notag \\
&&+\frac{(\Delta t)^{\alpha }}{2\Gamma (\alpha +3)}\sum_{i=2}^{m}\left[
f\left( t_{i},y_{i}\right) -2f\left( t_{i-1},y_{i-1}\right) +f\left(
t_{i-2},y_{i-2}\right) \right]  \notag \\
&&\times \left[ 
\begin{array}{c}
(m-i+1)^{\alpha }\left[ 
\begin{array}{c}
2(m-i)^{2}+(3\alpha +10)(m-i) \\ 
+2\alpha ^{2}+9\alpha +12%
\end{array}%
\right] \\ 
-(m-i)^{\alpha }\left[ 
\begin{array}{c}
2(m-i)^{2}+(5\alpha +10)(m-i) \\ 
+6\alpha ^{2}+18\alpha +12%
\end{array}%
\right]%
\end{array}%
\right] .  \label{3.25}
\end{eqnarray}%
In each iteration, the predictor (\ref{3.25}) is calculated and then
corrected by means of the implicit formula%
\begin{eqnarray}
y_{m+1} &=&y_{0}+\Upsilon _{m-1}+\frac{(\Delta t)^{\alpha }}{\Gamma (\alpha
+1)}f(t_{0},y_{0})\left[ (m+1)^{\alpha }-m^{\alpha }\right]  \notag \\
&&+\frac{(\Delta t)^{\alpha }}{\Gamma (\alpha +2)}\left(
f(t_{1},y_{1})-f(t_{0},y_{0})\right) \left[ (m+1)^{\alpha +1}-m^{\alpha
+1}-(\alpha +1)m^{\alpha }\right]  \notag \\
&&+\frac{(\Delta t)^{\alpha }}{\Gamma (\alpha +1)}f\left(
t_{m+1},y_{m+1}^{P}\right) +\frac{\alpha (\Delta t)^{\alpha }}{\Gamma
(\alpha +2)}\left( f\left( t_{m},y_{m}\right) -f\left(
t_{m+1},y_{m+1}^{P}\right) \right)  \notag \\
&&-\frac{\alpha (\Delta t)^{\alpha }}{2\Gamma (\alpha +3)}\left( f\left(
t_{m+1},y_{m+1}^{P}\right) -2f\left( t_{m},y_{m}\right) +f\left(
t_{m-1},y_{m-1}\right) \right) .  \label{3.26}
\end{eqnarray}

\subsection{Caputo-Fabrizio Fractional Derivative}

In this section, we will follow the same steps to derive a
predictor-corrector numertical scheme for the Caputo-Fabrizio fractional
initial-value problem%
\begin{equation}
\left \{ 
\begin{array}{l}
_{0}^{CF}D_{t}^{\alpha }y(t)=f(t,y(t)),\medskip \\ 
y(0)=y_{0},%
\end{array}%
\right.  \label{4.1}
\end{equation}%
where the fractional order $\alpha \in (0,1)$ and $f$ is a nonlinear smooth
function chosen such that system (\ref{4.1}) admits a unique solution $y(t)$%
. Similar to the previous section, we start with the difference formula 
\begin{equation*}
y(t)-y(0)=\frac{1-\alpha }{M(\alpha )}f(t,y(t))+\frac{\alpha }{M(\alpha )}%
\int_{0}^{t}f(s,y(s))ds,
\end{equation*}%
which when evaluated at two points in time $t_{m}=m\Delta t$ and $%
t_{m+1}=(m+1)\Delta t$ yields%
\begin{equation*}
y\left( t_{m}\right) -y(0)=\frac{1-\alpha }{M(\alpha )}f\left( t_{m},y\left(
t_{m}\right) \right) +\frac{\alpha }{M(\alpha )}\int_{0}^{t_{m}}f(s,y(s))ds,
\end{equation*}%
and%
\begin{equation}
y\left( t_{m+1}\right) -y(0)=\frac{1-\alpha }{M(\alpha )}f\left(
t_{m+1},y\left( t_{m+1}\right) \right) +\frac{\alpha }{M(\alpha )}%
\int_{0}^{t_{m+1}}f(s,y(s))ds,  \label{4.2}
\end{equation}%
respectively. Taking the difference of the two points produces%
\begin{equation}
y\left( t_{m+1}\right) -y\left( t_{m}\right) =\frac{1-\alpha }{M(\alpha )}%
\left[ f\left( t_{m+1},y\left( t_{m+1}\right) \right) -f\left( t_{m},y\left(
t_{m}\right) \right) \right] +\frac{\alpha }{M(\alpha )}%
\int_{t_{m}}^{t_{m+1}}f(s,y(s))ds.  \label{4.3}
\end{equation}%
Function $f(s,y(s))$ can be approximated over the sub-interval $%
[t_{m},t_{m+1}]$ by means of the same second order Newton polynomial (\ref%
{2.4}), which was employed in the classical derivative case. The result is%
\begin{eqnarray}
y_{m+1} &=&y_{m}+\frac{\alpha \Delta t}{M(\alpha )}f\left(
t_{m+1},y_{m+1}\right) +\frac{1-\alpha }{M(\alpha )}\left[ f\left(
t_{m+1},y\left( t_{m+1}\right) \right) -f\left( t_{m},y\left( t_{m}\right)
\right) \right]  \notag \\
&&+\frac{\alpha }{M(\alpha )}\left( \frac{f\left( t_{m+1},y_{m+1}\right)
-f\left( t_{m},y_{m}\right) }{\Delta t}\right) \int_{t_{m}}^{t_{m+1}}\left(
s-t_{m+1}\right) ds  \notag \\
&&+\frac{\alpha }{M(\alpha )}\left( \frac{f\left( t_{m+1},y_{m+1}\right)
-2f\left( t_{m},y_{m}\right) +f\left( t_{m-1},y_{m-1}\right) }{2(\Delta
t)^{2}}\right) \int_{t_{m}}^{t_{m+1}}\left( s-t_{m}\right) \left(
s-t_{m+1}\right) ds.  \label{4.4}
\end{eqnarray}%
Replacing the integrals by their respective values from (\ref{2.5.1}) and (%
\ref{2.5.2}) leads to the formula%
\begin{eqnarray}
y_{m+1} &=&y_{m}+\frac{1-\alpha }{M(\alpha )}\left[ f\left( t_{m+1},y\left(
t_{m+1}\right) \right) -f\left( t_{m},y\left( t_{m}\right) \right) \right] 
\notag \\
&&+\frac{\alpha \Delta t}{M(\alpha )}f\left( t_{m+1},y_{m+1}\right) -\left[
f\left( t_{m+1},y_{m+1}\right) -f\left( t_{m},y_{m}\right) \right] \frac{%
\alpha \Delta t}{2M(\alpha )}  \notag \\
&&-\left[ f\left( t_{m+1},y_{m+1}\right) -2f\left( t_{m},y_{m}\right)
+f\left( t_{m-1},y_{m-1}\right) \right] \frac{\alpha \Delta t}{12M(\alpha )}.
\label{4.5}
\end{eqnarray}%
Again, the terms $y_{m+1}$ appearing on the right hand side of the implicit
formula (\ref{4.5}) are replaced by the prediction $y_{m+1}^{P}$ obtained
using the Atangana-Seda scheme developed for the Caputo-Fabrizio fractional
derivative in \cite{Atangana2020,Atangana2020Corrigendum}. This yields the implicit corrector formula%
\begin{eqnarray}
y_{m+1}= &&y_{m}+\frac{1-\alpha }{M(\alpha )}\left[ f\left(
t_{m+1},y_{m+1}^{P}\right) -f\left( t_{m},y_{m}\right) \right]  \notag \\
&&+\frac{\alpha }{M(\alpha )}\left[ \frac{5}{12}f\left(
t_{m+1},y_{m+1}^{P}\right) \Delta t+\frac{2}{3}f\left( t_{m},y_{m}\right)
\Delta t-f\left( t_{m-1},y_{m-1}\right) \frac{\Delta t}{12}\right] ,
\label{4.6}
\end{eqnarray}%
with the predictor%
\begin{eqnarray}
y_{m+1}^{P} &=&y_{m}+\frac{1-\alpha }{M(\alpha )}\left[ f\left(
t_{m},y_{m}\right) -f\left( t_{m-1},y_{m-1}\right) \right]  \notag \\
&&+\frac{\alpha }{M(\alpha )}\left[ \frac{5}{12}f\left(
t_{m-2},y_{m-2}\right) \Delta t-\frac{4}{3}f\left( t_{m-1},y_{m-1}\right)
\Delta t+\frac{23}{12}f\left( t_{m},y_{m}\right) \Delta t\right] .
\label{4.7}
\end{eqnarray}

\subsection{Atangana-Baleanu Fractional Derivative}

The third type of fractional derivative we would like to consider is the
Atangana-Baleanu derivative. Let us consider the initial-value problem%
\begin{equation}
\left \{ 
\begin{array}{l}
_{0}^{ABC}D_{t}^{\alpha }y(t)=f(t,y(t)),\medskip \\ 
y(0)=y_{0},%
\end{array}%
\right.  \label{5.1}
\end{equation}%
where, as usual, the fractional order $\alpha \in (0,1)$ and $f$ is a smooth
nonlinear function that guarantees the existence of a unique solution $y(t)$
for (\ref{5.1}). In order to obtain a predictor-corrector numerical scheme
that solves (\ref{5.1}), we use the Atangana-Baleanu integral to produce%
\begin{equation*}
y(t)-y(0)=\frac{1-\alpha }{AB(\alpha )}f(t,y(t))+\frac{\alpha }{AB(\alpha
)\Gamma (\alpha )}\int_{0}^{t}f(s,y(s))(t-s)^{\alpha -1}ds,
\end{equation*}%
which leads to the approximation of $y(t)$ at $t_{m+1}=(m+1)\Delta t$ given
by%
\begin{equation}
y\left( t_{m+1}\right) =y(0)+\frac{1-\alpha }{AB(\alpha )}f\left(
t_{m+1},y\left( t_{m+1}\right) \right) +\frac{\alpha }{AB(\alpha )\Gamma
(\alpha )}\sum_{i=0}^{m}\int_{t_{i}}^{t_{i+1}}f(s,y(s))\left(
t_{m+1}-s\right) ^{\alpha -1}ds,  \label{5.2}
\end{equation}%
where $t_{0}=0$. Using the Newton polynomial (\ref{3.4}) to approximate
function $f(s,y(s))$ in (\ref{5.2}) yields%
\begin{eqnarray}
y(t_{m+1}) &=&y(0)+\frac{1-\alpha }{AB(\alpha )}f\left( t_{m+1},y\left(
t_{m+1}\right) \right)  \notag \\
&&+\frac{\alpha }{AB(\alpha )\Gamma (\alpha )}\int_{0}^{t_{1}}\left[
f(t_{0},y(t_{0}))+\left( \frac{f(t_{1},y(t_{1}))-f(t_{0},y(t_{0}))}{\Delta t}%
\right) s\right] \left( t_{m+1}-s\right) ^{\alpha -1}ds  \notag \\
&&+\frac{\alpha }{AB(\alpha )\Gamma (\alpha )}\sum_{i=1}^{m}%
\int_{t_{i}}^{t_{i+1}}\left \{ 
\begin{array}{l}
f\left( t_{i+1},y\left( t_{i+1}\right) \right) \\ 
+\frac{f\left( t_{i+1},y\left( t_{i+1}\right) \right) -f\left( t_{i},y\left(
t_{i}\right) \right) }{\Delta t}\left( s-t_{i+1}\right) \\ 
+\frac{f\left( t_{i+1},y\left( t_{i+1}\right) \right) -2f\left(
t_{i},y\left( t_{i}\right) \right) +f\left( t_{i-1},y\left( t_{i-1}\right)
\right) }{2(\Delta t)^{2}} \\ 
\times \left( s-t_{i}\right) \left( s-t_{i+1}\right)%
\end{array}%
\right \} \left( t_{m+1}-s\right) ^{\alpha -1}ds,  \label{5.3}
\end{eqnarray}%
which can be simplified and rearranged to the form%
\begin{eqnarray}
y_{m+1} &=&y_{0}+\frac{1-\alpha }{AB(\alpha )}f\left( t_{m+1},y\left(
t_{m+1}\right) \right)  \notag \\
&&+\frac{\alpha }{AB(\alpha )\Gamma (\alpha )}f(t_{0},y_{0})\int_{0}^{t_{1}}%
\left( t_{m+1}-s\right) ^{\alpha -1}ds  \notag \\
&&+\frac{\alpha }{AB(\alpha )\Gamma (\alpha )}\left( \frac{%
f(t_{1},y_{1})-f(t_{0},y_{0})}{\Delta t}\right) \int_{0}^{t_{1}}s\left(
t_{m+1}-s\right) ^{\alpha -1}ds  \notag \\
&&+\frac{\alpha }{AB(\alpha )\Gamma (\alpha )}\sum_{i=1}^{m}f\left(
t_{i+1},y_{i+1}\right) \int_{t_{i}}^{t_{i+1}}\left( t_{m+1}-s\right)
^{\alpha -1}ds  \notag \\
&&+\frac{\alpha }{AB(\alpha )\Gamma (\alpha )}\sum_{i=1}^{m}\frac{f\left(
t_{i+1},y_{i+1}\right) -f\left( t_{i},y_{i}\right) }{\Delta t}%
\int_{t_{i}}^{t_{i+1}}\left( s-t_{i+1}\right) \left( t_{m+1}-s\right)
^{\alpha -1}ds  \notag \\
&&+\frac{\alpha }{AB(\alpha )\Gamma (\alpha )}\sum_{i=1}^{m}\frac{f\left(
t_{i+1},y_{i+1}\right) -2f\left( t_{i},y_{i}\right) +f\left(
t_{i-1},y_{i-1}\right) }{2(\Delta t)^{2}}  \notag \\
&&\times \int_{t_{i}}^{t_{i+1}}\left( s-t_{i}\right) \left( s-t_{i+1}\right)
\left( t_{m+1}-s\right) ^{\alpha -1}ds.  \label{5.4}
\end{eqnarray}%
Replacing the integrals with their respective values from (\ref{3.9})-(\ref%
{3.12}) leads to%
\begin{eqnarray}
y_{m+1} &=&y_{0}+\frac{1-\alpha }{AB(\alpha )}f\left( t_{m+1},y\left(
t_{m+1}\right) \right)  \notag \\
&&+\frac{\alpha (\Delta t)^{\alpha }}{AB(\alpha )\Gamma (\alpha +1)}%
f(t_{0},y_{0})\left[ (m+1)^{\alpha }-m^{\alpha }\right]  \notag \\
&&+\frac{\alpha (\Delta t)^{\alpha }}{AB(\alpha )\Gamma (\alpha +2)}\left(
f(t_{1},y_{1})-f(t_{0},y_{0})\right) \left[ (m+1)^{\alpha +1}-m^{\alpha
+1}-(\alpha +1)m^{\alpha }\right]  \notag \\
&&+\frac{\alpha (\Delta t)^{\alpha }}{AB(\alpha )\Gamma (\alpha +1)}%
\sum_{i=1}^{m}f\left( t_{i+1},y_{i+1}\right) \left[ (m-i+1)^{\alpha
}-(m-i)^{\alpha }\right]  \notag \\
&&+\frac{(\Delta t)^{\alpha }}{\Gamma (\alpha +2)}\sum_{i=1}^{m}\left(
f\left( t_{i+1},y_{i+1}\right) -f\left( t_{i},y_{i}\right) \right) \left[
(m-i-\alpha )(m-i+1)^{\alpha }-(m-i)^{\alpha +1}\right]  \notag \\
&&+\frac{(\Delta t)^{\alpha }}{2\Gamma (\alpha +3)}\sum_{i=1}^{m}\left(
f\left( t_{i+1},y_{i+1}\right) -2f\left( t_{i},y_{i}\right) +f\left(
t_{i-1},y_{i-1}\right) \right)  \notag \\
&&\times \left[ 
\begin{array}{l}
(m-i+1)^{\alpha }\left[ 2(m-i)^{2}-\alpha (m-i+1)+2(m-i)\right] \\ 
-(m-i)^{\alpha }\left[ 2(m-i)^{2}+\alpha (m-i)+2(m-i)\right]%
\end{array}%
\right] .  \label{5.5}
\end{eqnarray}%
Using the notation $\Upsilon _{m-1}$ defined earlier in (\ref{3.14})-(\ref%
{3.15}) and replacing the terms $y_{m+1}$ on the right hand side of the
formula by the predicted value $y_{m+1}^{P}$, we obtain the predictor-corrector method described by the implicit formula%
\begin{eqnarray}
y_{m+1} &=&y_{0}+\frac{1-\alpha }{AB(\alpha )}f\left(
t_{m+1},y_{m+1}^{P}\right) +\frac{\alpha }{AB(\alpha )}\Upsilon _{m-1}+\frac{%
\alpha (\Delta t)^{\alpha }}{AB(\alpha )\Gamma (\alpha +1)}f(t_{0},y_{0})%
\left[ (m+1)^{\alpha }-m^{\alpha }\right]  \notag \\
&&+\frac{\alpha (\Delta t)^{\alpha }}{AB(\alpha )\Gamma (\alpha +2)}\left(
f(t_{1},y_{1})-f(t_{0},y_{0})\right) \left[ (m+1)^{\alpha +1}-m^{\alpha
+1}-(\alpha +1)m^{\alpha }\right]  \notag \\
&&+\frac{\alpha (\Delta t)^{\alpha }}{AB(\alpha )\Gamma (\alpha +1)}f\left(
t_{m+1},y_{m+1}^{P}\right) +\frac{\alpha ^{2}(\Delta t)^{\alpha }}{AB(\alpha
)\Gamma (\alpha +2)}\left( f\left( t_{m},y_{m}\right) -f\left(
t_{m+1},y_{m+1}^{P}\right) \right)  \notag \\
&&-\frac{\alpha ^{2}(\Delta t)^{\alpha }}{2AB(\alpha )\Gamma (\alpha +3)}%
\left( f\left( t_{m+1},y_{m+1}^{P}\right) -2f\left( t_{m},y_{m}\right)
+f\left( t_{m-1},y_{m-1}\right) \right) ,  \label{5.6}
\end{eqnarray}%
with the improved explicit Atangana-Seda predictor%
\begin{eqnarray}
y_{m+1}^{P} &=&y_{0}+\frac{1-\alpha }{AB(\alpha )}f\left( t_{m},y_{m}\right)
+\frac{\alpha (\Delta t)^{\alpha }}{AB(\alpha )\Gamma (\alpha +1)}%
\sum_{i=0}^{1}f\left( t_{i},y_{i}\right) \left[ (m-i+1)^{\alpha
}-(m-i)^{\alpha }\right]  \notag \\
&&+\frac{\alpha (\Delta t)^{\alpha }}{AB(\alpha )\Gamma (\alpha +2)}%
\sum_{i=0}^{1}\left( f\left( t_{i+1},y_{i+1}\right) -f\left(
t_{i},y_{i}\right) \right)  \notag \\
&&\times \left[ (m-i+1)^{\alpha +1}-(m-i)^{\alpha +1}-(\alpha
+1)(m-i)^{\alpha }\right]  \notag \\
&&+\frac{\alpha (\Delta t)^{\alpha }}{AB(\alpha )\Gamma (\alpha +1)}%
\sum_{i=2}^{m}f\left( t_{i-2},y_{i-2}\right) \left[ (m-i+1)^{\alpha
}-(m-i)^{\alpha }\right]  \notag \\
&&+\frac{\alpha (\Delta t)^{\alpha }}{AB(\alpha )\Gamma (\alpha +2)}%
\sum_{i=2}^{m}\left( f\left( t_{i-1},y_{i-1}\right) -f\left(
t_{i-2},y_{i-2}\right) \right)  \notag \\
&&\times \left[ 
\begin{array}{c}
(m-i+1)^{\alpha }(m-i+3+2\alpha ) \\ 
-(m-i)^{\alpha }(m-i+3+3\alpha )%
\end{array}%
\right]  \notag \\
&&+\frac{\alpha (\Delta t)^{\alpha }}{2AB(\alpha )\Gamma (\alpha +3)}%
\sum_{i=2}^{m}\left[ f\left( t_{i},y_{i}\right) -2f\left(
t_{i-1},y_{i-1}\right) +f\left( t_{i-2},y_{i-2}\right) \right]  \notag \\
&&\times \left[ 
\begin{array}{l}
(m-i+1)^{\alpha }\left[ 2(m-i)^{2}+(3\alpha +10)(m-i)+2\alpha ^{2}+9\alpha
+12\right] \\ 
-(m-i)^{\alpha }\left[ 2(m-i)^{2}+(5\alpha +10)(m-i)+6\alpha ^{2}+18\alpha
+12\right]%
\end{array}%
\right] .  \label{5.7}
\end{eqnarray}%
Note that this predictor is obtained in the same was as that of the Caputo
derivative in Section \ref{Section_Caputo}.

\subsection{Concluding Remarks}

\begin{remark}
The predictor term $y_{m+1}^{P}$ used in each\ of the previous scenarios\
can be replaced by any other scheme including, for instance, the ones in 
\cite{Diethelm1998,Odibat2008}. In the cases of the
Caputo-Fabrizio/Atangana-Baleanu fractional derivatives, some minor
modifications would have to be made to the methods.
\end{remark}

\begin{remark}
In the initial-value problem (\ref{3.1}), if $n-1<\alpha \leqslant n\in 
\mathbb{N}$ and $y_{0}=\left( y_{0,1},\dots ,y_{0,n}\right) $, then the
initial value $y_{0}$ on the right hand side of (\ref{3.25}) and (\ref{3.26}%
) needs to be replaced by the sum%
\begin{equation*}
\sum_{k=0}^{n-1}\frac{t_{m+1}^{k}}{k!}y_{0,k+1}.
\end{equation*}
\end{remark}

\section{Numerical Experiments}

In this section, we will present simulation results obtained by means of the
predictor-corrector numerical methods proposed in this paper for different
initial value problems. In the last example, we will consider a fractional
activator-inhibitor Gierer-Meinhardt model whose dynamics are to be analyzed
based on the obtained numerical solutions.

\begin{example}
We start with the classical initial-value problem%
\begin{equation}
\left \{ 
\begin{array}{l}
\frac{dy(t)}{dt}=2y(t)+3,\text{\medskip } \\ 
y(0)=1,%
\end{array}%
\right.  \label{6.1}
\end{equation}%
which has the exact solution%
\begin{equation}
y(t)=\frac{5}{2}e^{2t}-\frac{3}{2}.  \label{6.2}
\end{equation}%
Figure \ref{Fig1_Example1} depicts the exact solution (\ref{6.2}) along with
the numerical solutions obtained by means of the proposed method and the
standard Atangana-Seda method. The absolute error results are shown in Table %
\ref{Example1_errTable} for different values of the numerical step size. We
see that the proposed method for the classical derivative given in (\ref{2.7}%
) as well as the Caputo method in (\ref{3.26}) applied with $\alpha =1$
achieve a considerably lower error than the Atangana-Seda and two-step
Adams-Bashforth methods. 
\begin{table}[th]
\caption{Comparison of the maximum absolute errors of various numerical
methods for problem (\protect \ref{6.1}) with $t\in \lbrack 0,1]$.}
\label{Example1_errTable}\centering
\begin{tabular}{lcccc}
\toprule Method & $\Delta t=\frac{1}{16}$ & $\Delta t=\frac{1}{64}$ & $%
\Delta t=\frac{1}{200}$ & $\Delta t=\frac{1}{1024}$ \\ 
\midrule Proposed PC (\ref{3.26}) fractional, $\alpha =1$ & $2.6019\times
10^{-3}$ & $7.8442\times 10^{-5}$ & $2.9104\times 10^{-6}$ & $2.2690\times
10^{-8}$ \\ 
Proposed PC (\ref{2.7}) & $4.7391\times 10^{-3}$ & $9.6052\times 10^{-5}$ & $%
3.3246\times 10^{-6}$ & $2.5281\times 10^{-8}$ \\ 
Atangana-Seda \cite{Atangana2020} & $2.0657\times 10^{-2}$ & $3.9611\times
10^{-4}$ & $1.3570\times 10^{-5}$ & $1.0281\times 10^{-7}$ \\ 
Two-step Adams-Bashforth & $2.0503\times 10^{-1}$ & $1.4503\times 10^{-2}$ & 
$1.5223\times 10^{-3}$ & $5.8597\times 10^{-5}$ \\ 
\bottomrule &  &  &  & 
\end{tabular}%
\end{table}

\begin{figure}[tbp]
\centering
\includegraphics[width=11cm]{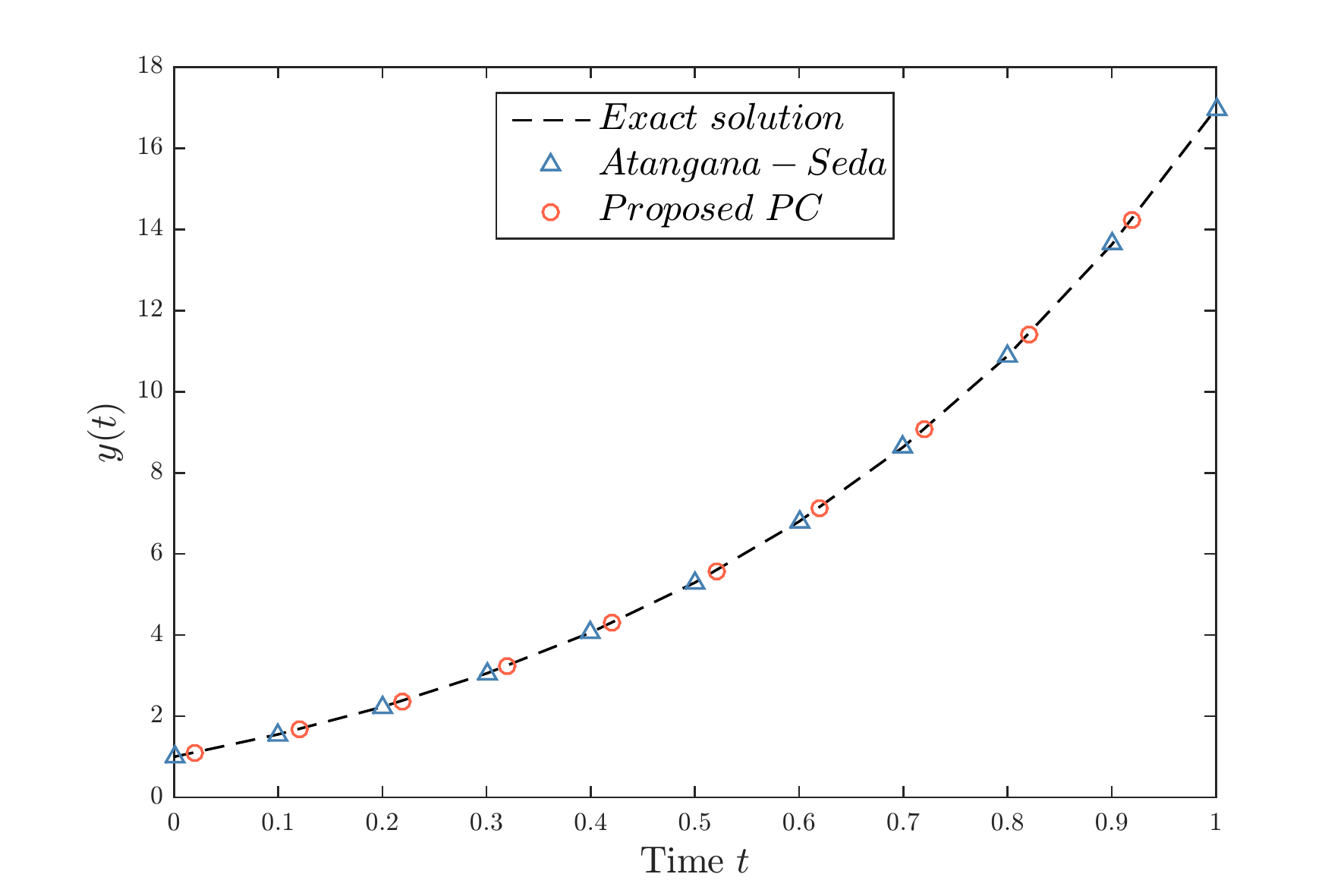}
\caption{Solution of problem (\protect \ref{6.1}) for $t\in \lbrack 0,1]$.}
\label{Fig1_Example1}
\end{figure}
\end{example}

\begin{example}
Let us consider another initial-value problem with a classical derivative:%
\begin{equation}
\left \{ 
\begin{array}{l}
\frac{dy(t)}{dt}=-\cos (2t)y^{2}(t),\text{\medskip } \\ 
y(0)=1.%
\end{array}%
\right.  \label{6.3}
\end{equation}%
The exact solution of this problem is known to be%
\begin{equation}
y(t)=\frac{2}{2+\sin (2t)}.  \label{6.4}
\end{equation}%
The exact solution (\ref{6.4}) is depicted in Figure \ref{Fig2_Example2}
alongside the numerical solution obtained by means of the proposed numerical
scheme (\ref{2.7}) and the Atangana-Seda solution. The error performance is
detailed in\ Table \ref{Example2_errTable}. Again, the proposed schemes
achieve a noticeably superior performance. 
\begin{table}[th]
\caption{Comparison of the maximum absolute errors of various numerical
methods for problem (\protect \ref{6.3}) with $t\in \lbrack 0,30]$. }
\label{Example2_errTable}\centering
\begin{tabular}{lcccc}
\toprule Method & $\Delta t=\frac{1}{16}$ & $\Delta t=\frac{1}{64}$ & $%
\Delta t=\frac{1}{200}$ & $\Delta t=\frac{1}{700}$ \\ 
\midrule Proposed PC (\ref{3.26}) fractional, $\alpha =1$ & $8.3152\times
10^{-3}$ & $2.2772\times 10^{-5}$ & $6.5114\times 10^{-7}$ & $2.6151\times
10^{-8}$ \\ 
Proposed PC (\ref{2.7}) & $8.9834\times 10^{-3}$ & $1.0474\times 10^{-4}$ & $%
3.1725\times 10^{-6}$ & $7.1930\times 10^{-8}$ \\ 
Atangana-Seda \cite{Atangana2020} & $2.2712\times 10^{-2}$ & $3.4369\times
10^{-4}$ & $1.1236\times 10^{-5}$ & $2.6193\times 10^{-7}$ \\ 
Two-step Adams-Bashforth & $2.1387\times 10^{-2}$ & $1.3589\times 10^{-3}$ & 
$1.3984\times 10^{-4}$ & $1.1436\times 10^{-5}$ \\ 
\bottomrule &  &  &  & 
\end{tabular}%
\end{table}

\begin{figure}[tbp]
\centering
\includegraphics[width=11cm]{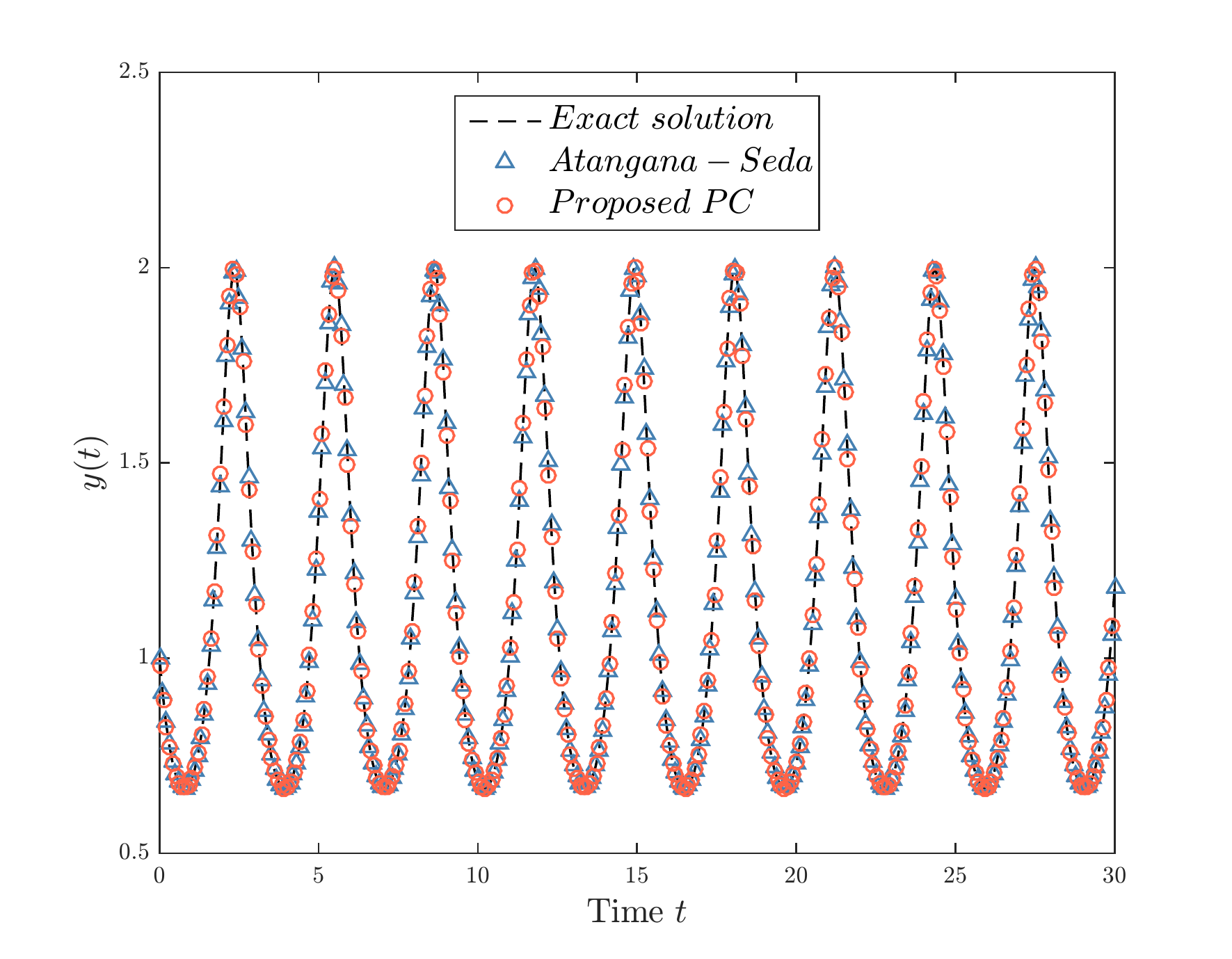}
\caption{Solution of problem (\protect \ref{6.3}) for $t\in \lbrack 0,30]$.}
\label{Fig2_Example2}
\end{figure}
\end{example}

\begin{example}
Next, we consider the fractional Caputo initial-value problem%
\begin{equation}
\left \{ 
\begin{array}{l}
\ _{0}^{C}D_{t}^{\alpha }y(t)=t^{\beta },\text{\medskip } \\ 
\ y(0)=0,%
\end{array}%
\right.  \label{6.5}
\end{equation}%
for some real constant $\beta $,\ which admits the unique exact solution%
\begin{equation}
y(t)=\frac{\Gamma (\beta +1)}{\Gamma (\alpha +\beta +1)}t^{\alpha +\beta }.
\label{6.6}
\end{equation}%
Figure \ref{Fig3_Example3} shows the exact solution (\ref{6.6}) along with
the numerical solution obtained by means of the proposed predictor corrector
scheme (\ref{3.26}) and the standard and improved Atangana-Seda methods for $%
\beta =0.9$ and $\alpha \in \left \{ 0.25,0.87\right \} $. The absolute
error results are presented in Table \ref{Example3_errTable} for the same
value of $\beta $ and $\alpha \in \left \{ 0.25,0.56,0.87\right \} $ with
different numerical step sizes. In all scenratios, the absolute error
achieved by the proposed method is lower than the improved Atangana-Seda
method, which in turn is lower than the standard one. 
\begin{table}[th]
\caption{Comparison of the maximum absolute errors of various methods for
problem (\protect \ref{6.5}) with $\protect \beta =0.9$ and $t\in \lbrack 0,3]$%
.}
\label{Example3_errTable}\centering
\begin{tabular}{lcccccc}
\toprule & \multicolumn{2}{c}{$\alpha =0.25$} & \multicolumn{2}{c}{$\alpha
=0.56$} & \multicolumn{2}{c}{$\alpha =0.87$} \\ 
\cmidrule(lr){2-3} \cmidrule(lr){4-5} \cmidrule(lr){6-7} Method & $\Delta t=%
\frac{1}{100}$ & $\Delta t=\frac{1}{800}$ & $\Delta t=\frac{1}{100}$ & $%
\Delta t=\frac{1}{400}$ & $\Delta t=\frac{1}{100}$ & $\Delta t=\frac{1}{200}$
\\ 
\midrule PPC (\ref{3.26}) & $6.8792\times 10^{-5}$ & $6.2948\times 10^{-6}$
& $2.8000\times 10^{-5}$ & $3.6996\times 10^{-6}$ & $7.6132\times 10^{-6}$ & 
$4.4095\times 10^{-7}$ \\ 
IAS (\ref{3.25}) & $3.9492\times 10^{-4}$ & $3.6137\times 10^{-5}$ & $%
8.4439\times 10^{-5}$ & $1.1157\times 10^{-5}$ & $1.9429\times 10^{-5}$ & $%
1.1253\times 10^{-6}$ \\ 
AS \cite{Atangana2020} & $2.3016\times 10^{-3}$ & $2.1060\times 10^{-4}$ & $%
1.2783\times 10^{-3}$ & $1.6891\times 10^{-4}$ & $4.9855\times 10^{-4}$ & $%
2.8876\times 10^{-5}$ \\ 
\bottomrule &  &  &  &  &  & 
\end{tabular}%
\end{table}

\begin{figure}[tbp]
\centering
\includegraphics[width=14cm]{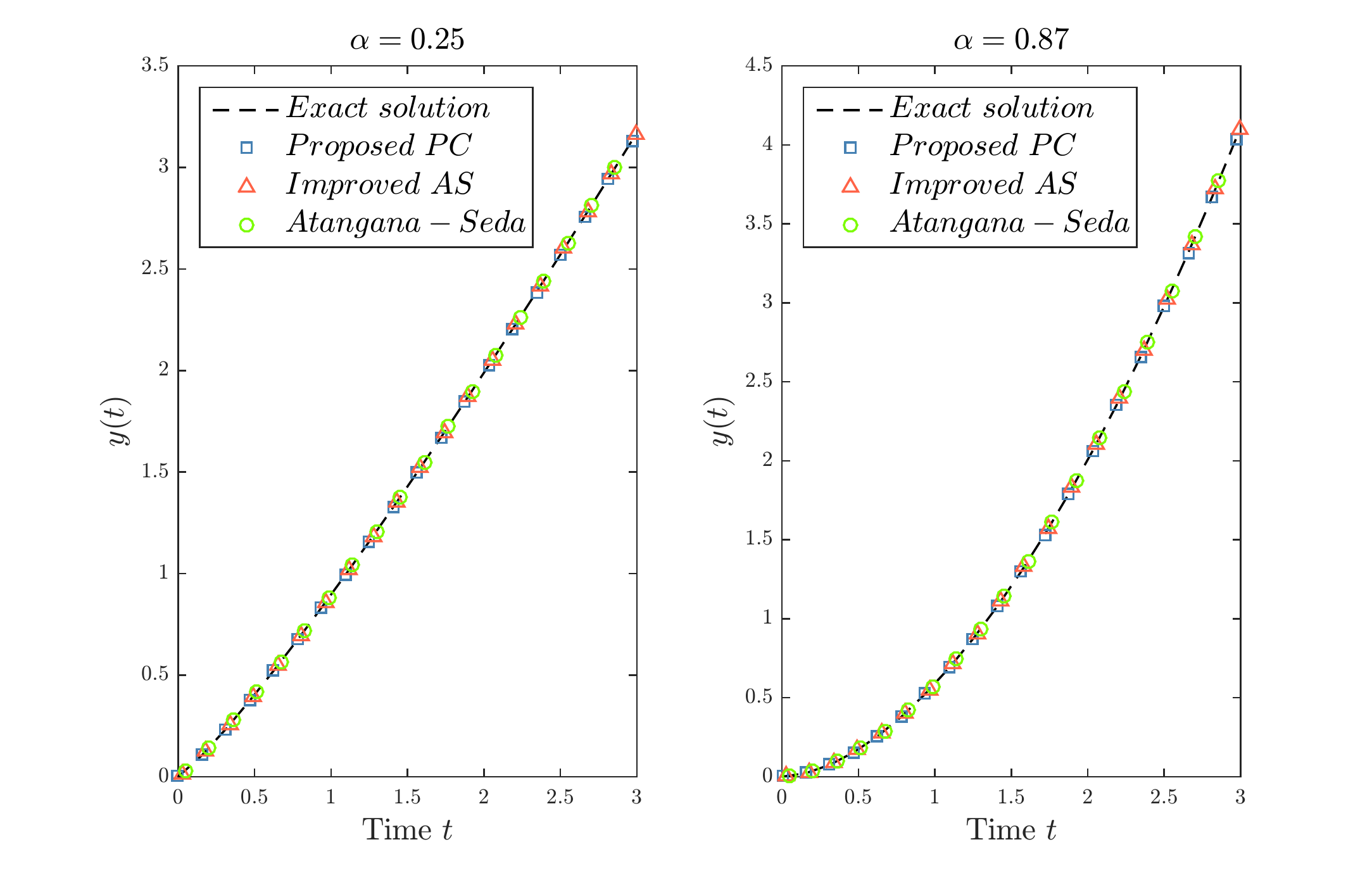}
\caption{Solution of problem (\protect \ref{6.5}) for $\protect \beta =0.9$
and $t\in \lbrack 0,3]$.}
\label{Fig3_Example3}
\end{figure}
\end{example}

\begin{example}
Let us consider the\ fractional Caputo initial-value problem%
\begin{equation}
\left \{ 
\begin{array}{l}
\ _{0}^{C}D_{t}^{\alpha }y(t)=\frac{2t^{2-\alpha }}{\Gamma (3-\alpha )}-%
\frac{t^{1-\alpha }}{\Gamma (2-\alpha )}-y(t)-t+t^{2},\text{\medskip } \\ 
\ y(0)=0.%
\end{array}%
\right.  \label{6.7}
\end{equation}%
The exact solution of (\ref{6.7}) can be shown to be%
\begin{equation}
y(t)=t^{2}-t.  \label{6.8}
\end{equation}%
Figure \ref{Fig4_Example4}\ and Table \ref{Example4_errTable}\ present the
numerical solutions of (\ref{6.7}) in comparison to the exact solution (\ref%
{6.8}) for different fractional orders and numerical steps sizes. Again, the
proposed method (\ref{3.26}) is superior to the Atangana-Seda method and the
improved method (\ref{3.25}). 
\begin{table}[th]
\caption{Comparison of the maximum absolute errors of various methods for
problem (\protect \ref{6.7}) with $t\in \lbrack 0,1]$.}
\label{Example4_errTable}\centering%
\begin{tabular}{lcccccc}
\toprule & \multicolumn{2}{c}{$\alpha =0.4$} & \multicolumn{2}{c}{$\alpha
=0.65$} & \multicolumn{2}{c}{$\alpha =0.9$} \\ 
\cmidrule(lr){2-3} \cmidrule(lr){4-5} \cmidrule(lr){6-7} Method & $\Delta t=%
\frac{1}{64}$ & $\Delta t=\frac{1}{512}$ & $\Delta t=\frac{1}{64}$ & $\Delta
t=\frac{1}{512}$ & $\Delta t=\frac{1}{64}$ & $\Delta t=\frac{1}{512}$ \\ 
\midrule PPC (\ref{3.26}) & $7.6806\times 10^{-4}$ & $6.4455\times 10^{-5}$
& $3.1549\times 10^{-3}$ & $4.5513\times 10^{-4}$ & $6.4490\times 10^{-3}$ & 
$8.2593\times 10^{-4}$ \\ 
IAS (\ref{3.25}) & $5.7442\times 10^{-3}$ & $7.0486\times 10^{-4}$ & $%
8.8970\times 10^{-3}$ & $1.1129\times 10^{-3}$ & $1.2365\times 10^{-2}$ & $%
1.5685\times 10^{-3}$ \\ 
AS \cite{Atangana2020} & $1.5787\times 10^{-2}$ & $2.0038\times 10^{-3}$ & $%
2.5879\times 10^{-2}$ & $3.2837\times 10^{-3}$ & $3.4087\times 10^{-2}$ & $%
4.3577\times 10^{-3}$ \\ 
\bottomrule &  &  &  &  &  & 
\end{tabular}%
\end{table}

\begin{figure}[tbp]
\centering
\includegraphics[width=14cm]{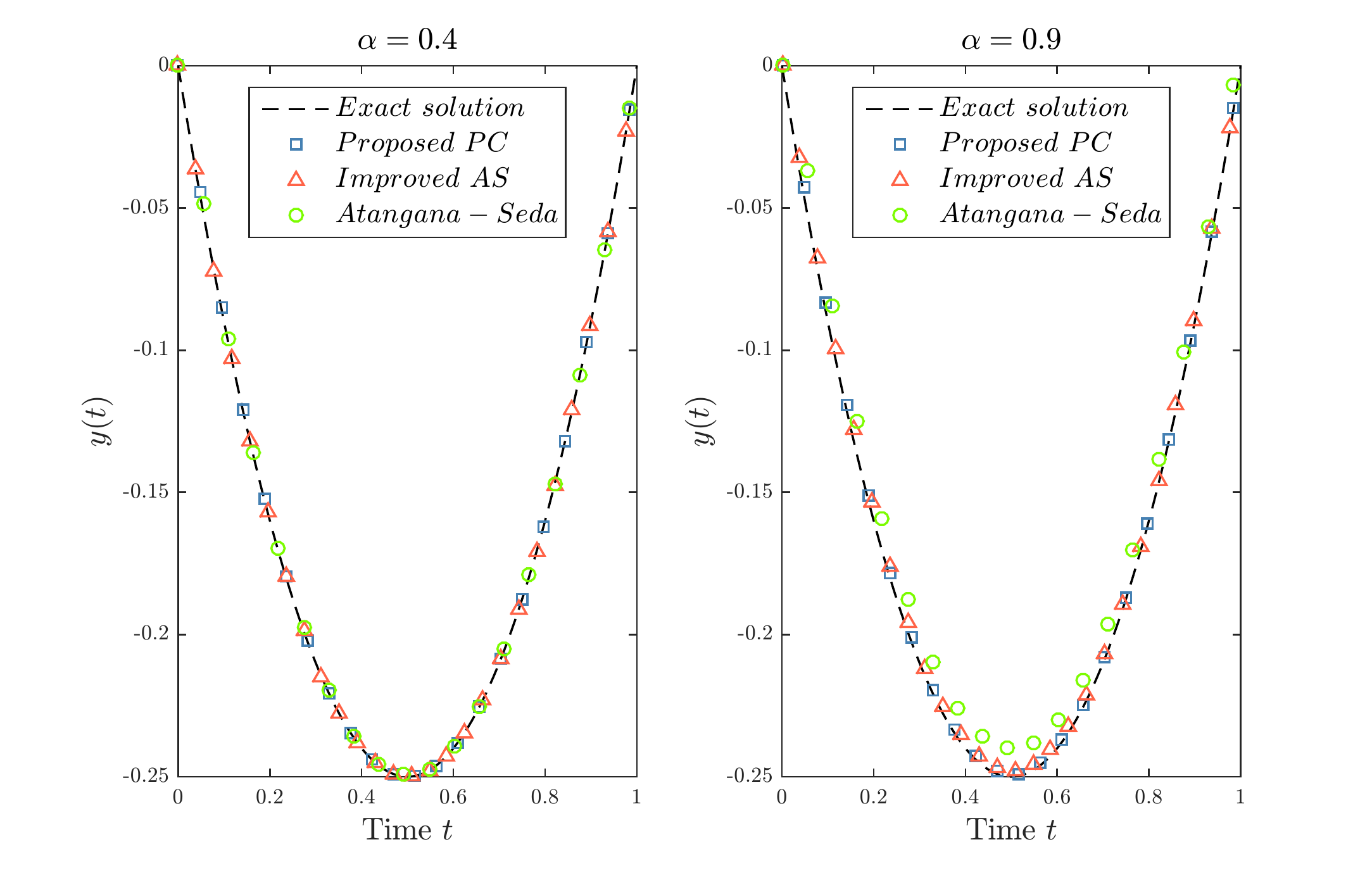}
\caption{Solution of problem (\protect \ref{6.7}) for $t\in \lbrack 0,1]$.}
\label{Fig4_Example4}
\end{figure}
\end{example}

\begin{example}
In the previous examples, we considered some simple single differential
equations with known exact solutions. Let us now analyze a realistic\
fractional activator-inhibitor model using analytical stability theory and
validate the theoretical results numerically by means of the proposed
method. Consider the system described by%
\begin{equation}
\left \{ 
\begin{array}{l}
_{0}^{C}D_{t}^{\alpha }a(t)=\varrho _{0}\varrho +c\varrho \frac{a(t)^{2}}{%
h(t)}-\mu a(t),\medskip  \\ 
_{0}^{C}D_{t}^{\alpha }h(t)=c^{\prime }\varrho ^{\prime }a(t)^{2}-\nu
h(t),\medskip  \\ 
a(0)=a_{0},\  \ h(0)=h_{0},%
\end{array}%
\right.   \label{6.9}
\end{equation}%
where $a(t)$ and $h(t)$ denote the concentrations of the activator and
inhibitor substances at time instant $t$, respectively. The constants $%
\varrho _{0}$,$\varrho $,$c$,$\mu $,$c^{\prime }$,$\varrho ^{\prime }$,$a_{0}
$,$h_{0}$ and $\nu $ are assumed to be positive real numbers, and the
fractional differentiation order $\alpha \in (0,1]$. For $\alpha =1$, system
(\ref{6.9}) reduces to the well known Gierer-Meinhardt model describing the
morphogenesis process \cite{Granero1984,Ruan1998}. Morphogenesis is the
biological process driving living organisms to take specific shapes.
Inclusion of a diffusion part in the Gierer-Meinhardt model was useful in
modeling the head formation of a fresh-water animal known as hydra \cite%
{GiererMeinhardt1972}. It is well established that system (\ref{6.9}) admits
the unique equilibrium point%
\begin{equation}
E^{\ast }=(a^{\ast },h^{\ast })  \label{6.10}
\end{equation}%
where 
\begin{equation}
a^{\ast }=\frac{\varrho _{0}\varrho c^{\prime }\varrho ^{\prime }+c\varrho
\nu }{\mu c^{\prime }\varrho ^{\prime }},  \label{6.11}
\end{equation}%
and%
\begin{equation}
h^{\ast }=\frac{c^{\prime }\varrho ^{\prime }}{\nu }(a^{\ast })^{2}.
\label{6.12}
\end{equation}%
Evaluating the\ Jacobian matrix of system (\ref{6.9}) at the unique
equilibrium $E^{\ast }$ yields%
\begin{equation}
J|_{E^{\ast }}=\left( 
\begin{array}{cc}
\frac{2c\mu \nu }{c\nu +c^{\prime }\varrho ^{\prime }\varrho _{0}}-\mu  & -%
\frac{c}{\varrho }\left( \frac{\mu \nu }{c\nu +c^{\prime }\varrho ^{\prime
}\varrho _{0}}\right) ^{2} \\ 
\frac{2\varrho \left( c\nu +c^{\prime }\varrho ^{\prime }\varrho _{0}\right) 
}{\mu } & -\nu 
\end{array}%
\right) .  \label{6.13}
\end{equation}%
The\ determinant and trace of the Jacobian are given by 
\begin{equation}
\mathrm{tr}J|_{E^{\ast }}=\frac{2\mu \nu c}{\nu c+\varrho _{0}\varrho
^{\prime }c^{\prime }}-\mu -\nu ,  \label{6.14}
\end{equation}%
and 
\begin{equation}
\mathrm{det}J|_{E^{\ast }}=\mu \nu ,  \label{6.15}
\end{equation}%
respectively. Hence, the characteristic equation of associated with $E^{\ast
}$ is%
\begin{equation}
\lambda ^{2}-\lambda \mathrm{tr}J|_{E^{\ast }}+\mathrm{det}J|_{E^{\ast }}=0,
\label{6.16}
\end{equation}%
leading\ to the eigenvalues%
\begin{equation}
\lambda _{1,2}=\frac{1}{2}\left( \mathrm{tr}J|_{E^{\ast }}\pm \sqrt{\mathrm{%
tr}^{2}J|_{E^{\ast }}-4\mathrm{det}J|_{E^{\ast }}}\right) .  \label{6.17}
\end{equation}%
The dynamics of (\ref{6.9}) can be analyzed by means of the results in \cite[%
Section 3]{Ahmed2007}. Firstly, if the discriminant of (\ref{6.16}) is equal
to zero, i.e. 
\begin{equation}
\mathrm{tr}^{2}J|_{E^{\ast }}-4\mathrm{det}J|_{E^{\ast }}=0,  \label{6.18}
\end{equation}%
the eigenvelues (\ref{6.17}) reduce to the real quantity%
\begin{equation}
\lambda _{1,2}=\frac{1}{2}\mathrm{tr}J|_{E^{\ast }}.  \label{6.19}
\end{equation}%
Hence, the equilibrium\ $E^{\ast }$ is asymptotically stable when $\mathrm{tr%
}J|_{E^{\ast }}<0$ and unstable when $\mathrm{tr}J|_{E^{\ast }}>0$ for all $%
\alpha \in (0,1]$.

Secondly, if the discriminant is strictly positive, i.e.%
\begin{equation}
\mathrm{tr}^{2}J|_{E^{\ast }}-4\mathrm{det}J|_{E^{\ast }}>0,  \label{6.20}
\end{equation}%
the eigenvalues (\ref{6.17}) are also real. However, we distinguish two
cases with respect to the asymptotic stability:

\begin{itemize}
\item If $\mathrm{tr}J|_{E^{\ast }}>0$, then%
\begin{equation}
\lambda _{1}=\frac{1}{2}\left( \mathrm{tr}J|_{E^{\ast }}+\sqrt{\mathrm{tr}%
^{2}J|_{E^{\ast }}-4\mathrm{det}J|_{E^{\ast }}}\right) >0.  \label{6.21}
\end{equation}%
Thus, $\left \vert \arg (\lambda _{1})\right \vert =0$ and $E^{\ast }$ is
unstable for all\ $\alpha \in (0,1]$.

\item If $\mathrm{tr}J|_{E^{\ast }}<0$, then%
\begin{equation}
\left \vert \arg (\lambda _{1,2})\right \vert =\pi >\frac{\alpha \pi }{2}\ 
\text{for }\alpha \in (0,1].  \label{6.22}
\end{equation}%
Thus, $E^{\ast }$ is asymptotically stable for all $\alpha \in (0,1]$.
\end{itemize}

Thirdly, if the discriminant is strictly negative, i.e.%
\begin{equation}
\mathrm{tr}^{2}J|_{E^{\ast }}-4\mathrm{det}J|_{E^{\ast }}<0,  \label{6.23}
\end{equation}%
the eigenvalues become%
\begin{equation}
\lambda _{1,2}=\frac{1}{2}\left( \mathrm{tr}J|_{E^{\ast }}\pm i\sqrt{4%
\mathrm{det}J|_{E^{\ast }}-\mathrm{tr}^{2}J|_{E^{\ast }}}\right) ,
\label{6.24}
\end{equation}%
leading to three distinguishable cases:

\begin{itemize}
\item If $\mathrm{tr}J|_{E^{\ast }}=0$, then%
\begin{equation}
\lambda _{1,2}=\pm i\sqrt{\mathrm{det}J|_{E^{\ast }}},  \label{6.25}
\end{equation}%
leading to%
\begin{equation}
\left \vert \arg (\lambda _{1,2})\right \vert =\frac{\pi }{2}>\frac{\alpha \pi 
}{2}\text{ for }\alpha \in (0,1).  \label{6.26}
\end{equation}%
Hence, $E^{\ast }$ is asymptotically stable for all $\alpha \in (0,1)$.

\item If $\mathrm{tr}J|_{E^{\ast }}<0$, then%
\begin{equation}
\left \vert \arg (\lambda _{1,2})\right \vert >\frac{\pi }{2}>\frac{\alpha \pi 
}{2}\text{ for }\alpha \in (0,1),  \label{6.27}
\end{equation}%
and, consequently, $E^{\ast }$ is asymptotically stable for all\ $\alpha \in
(0,1]$.

\item If $\mathrm{tr}J|_{E^{\ast }}>0$, then $E^{\ast }$ is asymptotically
stable for all $\alpha \in (0,1)$ if%
\begin{equation}
\tan ^{2}\left( \left \vert \mathrm{arg}\left( \lambda _{1,2}\right)
\right \vert \right) =\frac{4\mu \nu \left( c\nu +\varrho _{0}\varrho
^{\prime }c^{\prime }\right) ^{2}}{\left( c\nu (\mu -\nu )-\varrho
_{0}\varrho ^{\prime }c^{\prime }(\mu +\nu )\right) ^{2}}>\tan ^{2}\left( 
\frac{\alpha \pi }{2}\right) +1,  \label{6.28}
\end{equation}%
and unstable for all $\alpha \in (0,1)$ if%
\begin{equation}
\frac{4\mu \nu \left( c\nu +\varrho _{0}\varrho ^{\prime }c^{\prime }\right)
^{2}}{\left( c\nu (\mu -\nu )-\varrho _{0}\varrho ^{\prime }c^{\prime }(\mu
+\nu )\right) ^{2}}<\tan ^{2}\left( \frac{\alpha \pi }{2}\right) +1.
\label{6.29}
\end{equation}
\end{itemize}

\begin{remark}
\label{GM_remark_unstable}If the unique equilibrium $E^{\ast }$ of (\ref{6.9}%
) is unstable for some $\alpha \in (0,1)$, then $E^{\ast }$ is also unstable
for $\alpha =1$.
\end{remark}
Since an exact solution is not available for system (\ref%
{6.9}), visualizing the system dynamics requires numerical solutions, which
can be obtained using the proposed predictor-corrector method
described by (\ref{3.25})-(\ref{3.26}). The parameters adopted for the
simulations are listed in Table\  \ref{GM_table_par}. Condition (\ref{6.23})
can be easily verified and $\mathrm{tr}J|_{E^{\ast }}=\frac{6}{7}>0$. For $%
\alpha =0.85$, we have%
\begin{equation}
\tan ^{2}\left( \left \vert \mathrm{arg}\left( \lambda _{1,2}\right)
\right \vert \right) =\frac{392}{9}>\tan ^{2}\left( \frac{\alpha \pi }{2}%
\right) +1\approx 18.3497,
\end{equation}%
which implies that the equilibrium $E^{\ast }=(\frac{7}{4},\frac{49}{32})$
is asymptotically stable. The numerical solutions and corresponding phase
plot depicted in Figures \ref{GM_model_Fig1_stable} and \ref%
{GM_model_Fig2_phseP_stable}, respectively, agree with the theoretical
analysis as the solution converges towards $(\frac{7}{4},\frac{49}{32})$.
For $\alpha =0.95$, we have 
\begin{equation}
\tan ^{2}\left( \left \vert \mathrm{arg}\left( \lambda _{1,2}\right)
\right \vert \right) =\frac{392}{9}<\tan ^{2}\left( \frac{\alpha \pi }{2}%
\right) +1\approx 162.4476,
\end{equation}%
and thus, the equilibrium $E^{\ast }=(\frac{7}{4},\frac{49}{32})$ is
unstable. Again, the numerical results shown in in Figures \ref%
{GM_model_Fig3_Unstable} and \ref{GM_model_Fig4_phseP_stable} coincide with
the theoretical results as the solution is periodically stable around $(%
\frac{7}{4},\frac{49}{32})$. According to Remark \ref{GM_remark_unstable},
we conclude that the equilibrium $E^{\ast }$ of (\ref{6.9}) is unstable for $%
\alpha =1$. This result is confirmed by the numerical results depicted in
Figures \ref{GM_model_Fig5_1_Unstable} and \ref{GM_model_Fig6_2_phseP_stable}%
.
\begin{table}[th]
\caption{Parameter values of system (\protect \ref{6.9}) adopted in the
numerical simulations.}\centering
\begin{tabular}{ccccccccc}
\toprule$\varrho_{0}$ & $\varrho$ & $\mu$ & $\nu$ & $c$ & $\varrho^{\prime}$
& $c^{\prime}$ & $a_{0}$ & $h_{0}$ \\ 
\midrule $1$ & $1$ & $4$ & $2$ & $3 $ & $1$ & $1$ & $2$ & $3$ \\ 
\bottomrule &  &  &  &  &  &  &  & 
\end{tabular}%
\label{GM_table_par}
\end{table}

\begin{figure}[tbp]
\centering
\includegraphics[width=11cm]{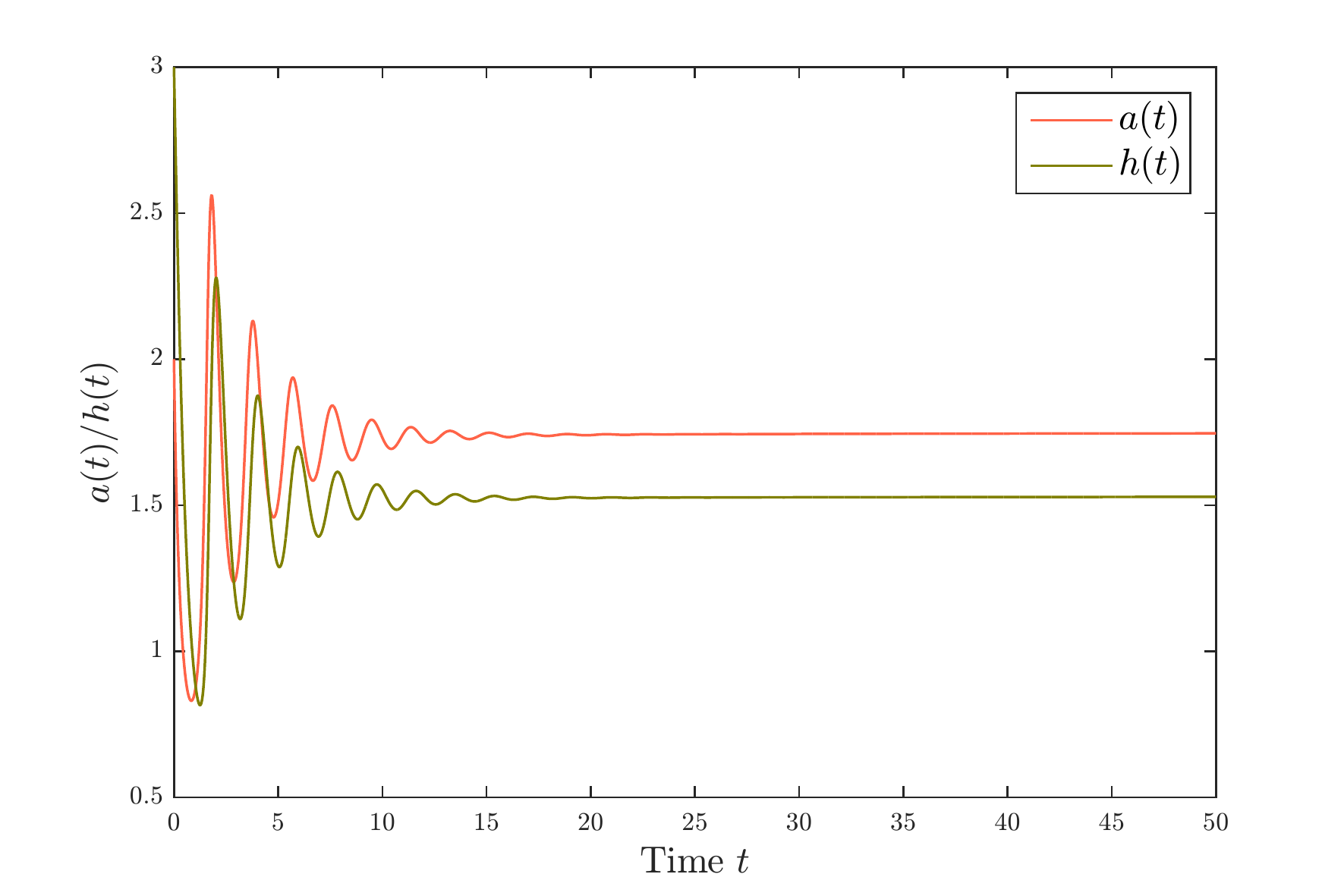}
\caption{The numerical solution of system (\protect \ref{6.9}) for $\protect%
\alpha =0.85$ with the parameters listed in Table \protect \ref{GM_table_par}.
}
\label{GM_model_Fig1_stable}
\end{figure}

\begin{figure}[tbp]
\centering
\includegraphics[width=11cm]{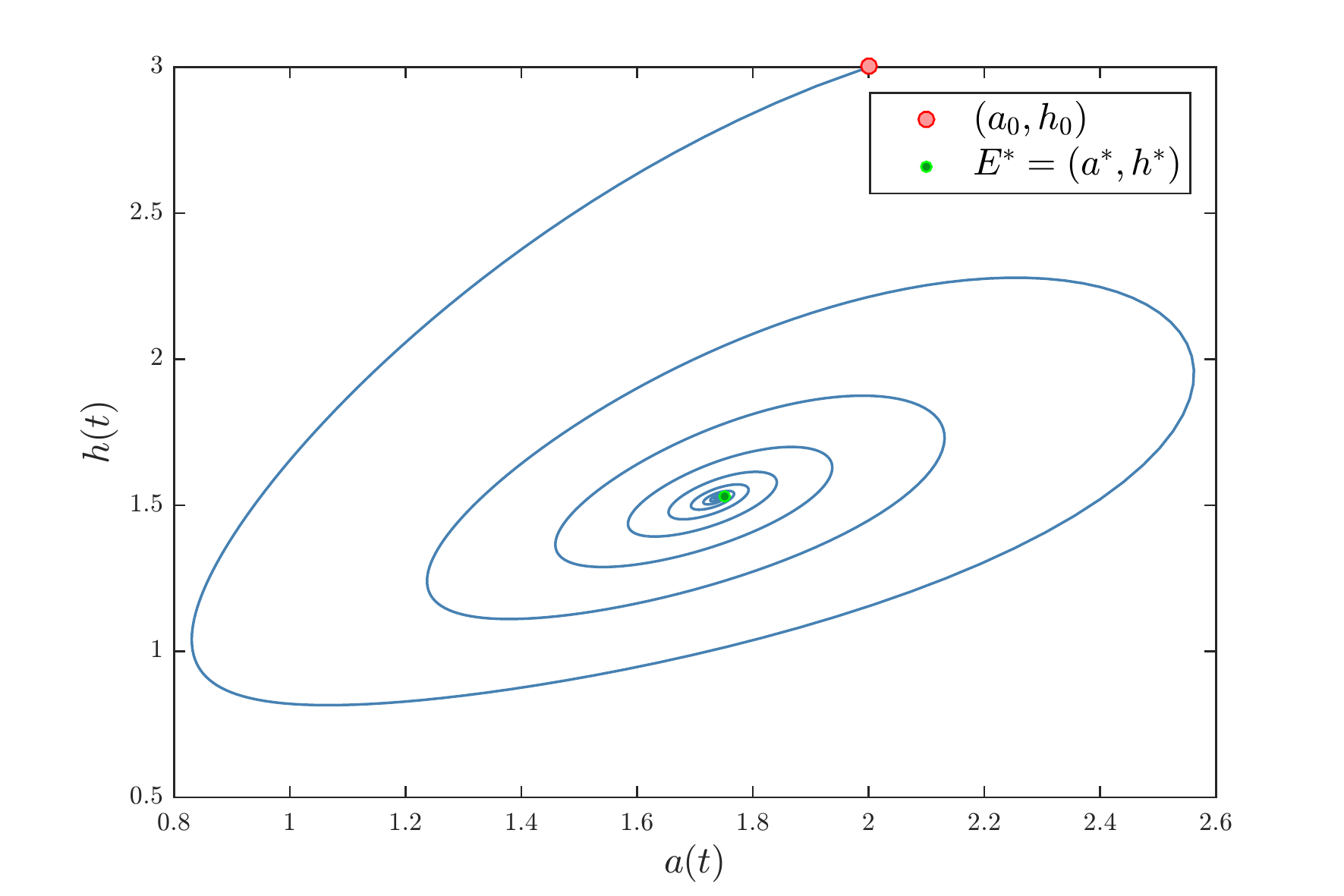}
\caption{Phase plot of system (\protect \ref{6.9}) for $\protect \alpha =0.85$
with the parameters listed in Table \protect \ref{GM_table_par}.}
\label{GM_model_Fig2_phseP_stable}
\end{figure}

\begin{figure}[tbp]
\centering
\includegraphics[width=11cm]{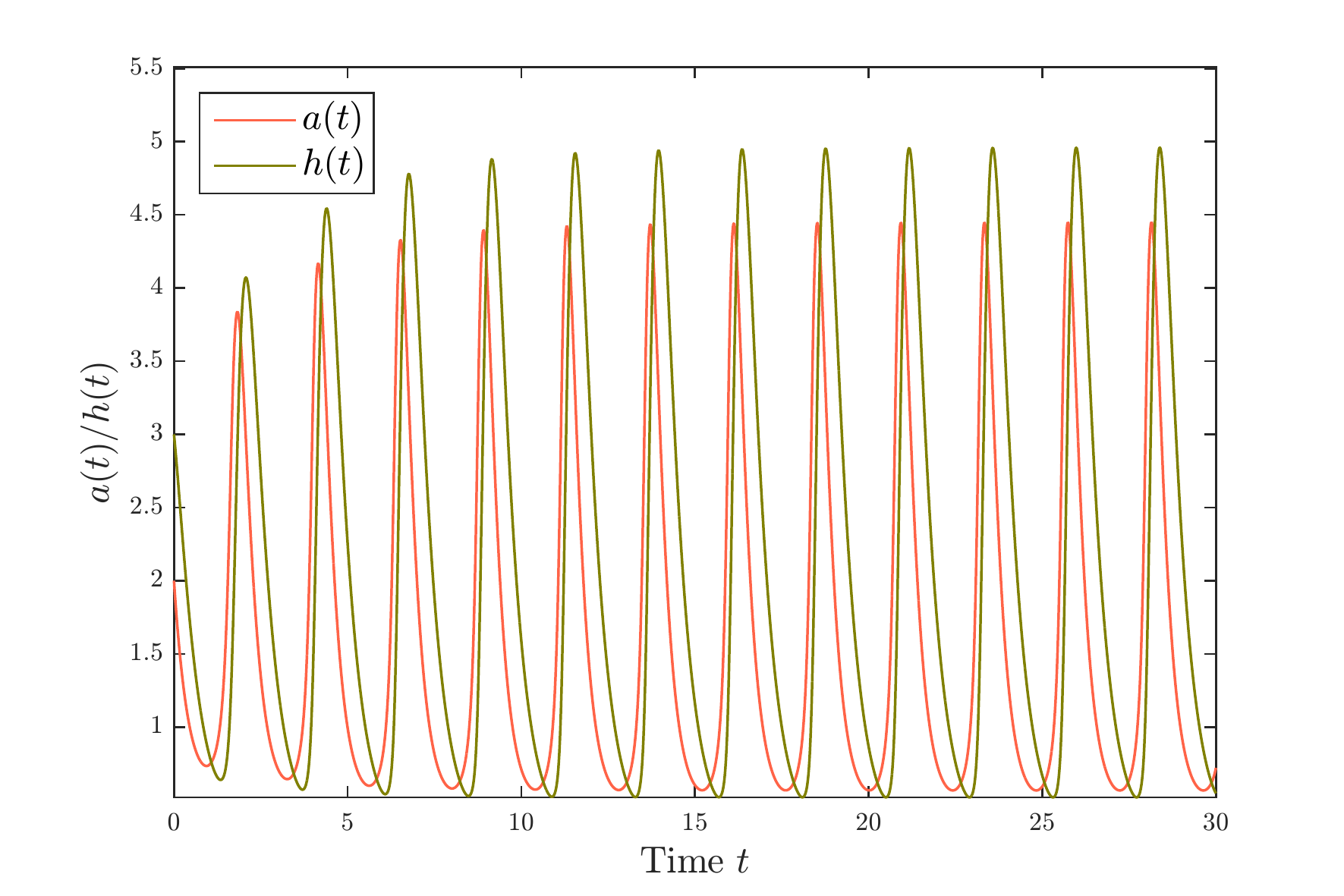}
\caption{The numerical solution of system (\protect \ref{6.9}) for $\protect%
\alpha =0.95$ with the parameters listed in Table \protect \ref{GM_table_par}.
}
\label{GM_model_Fig3_Unstable}
\end{figure}

\begin{figure}[tbp]
\centering
\includegraphics[width=11cm]{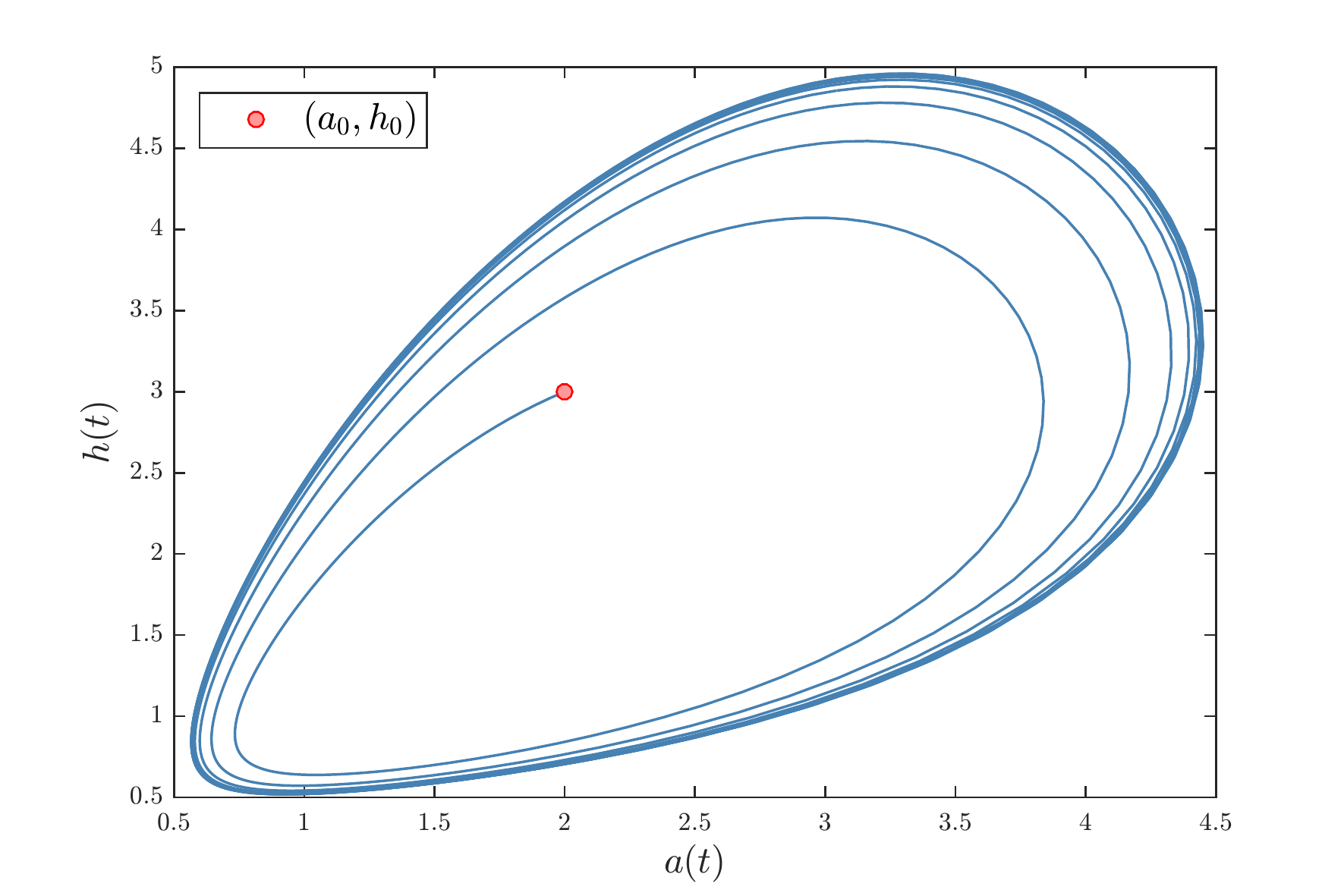}
\caption{Phase plot of system (\protect \ref{6.9}) for $\protect \alpha =0.95$
with the parameters listed in Table \protect \ref{GM_table_par}.}
\label{GM_model_Fig4_phseP_stable}
\end{figure}

\begin{figure}[tbp]
\centering
\includegraphics[width=11cm]{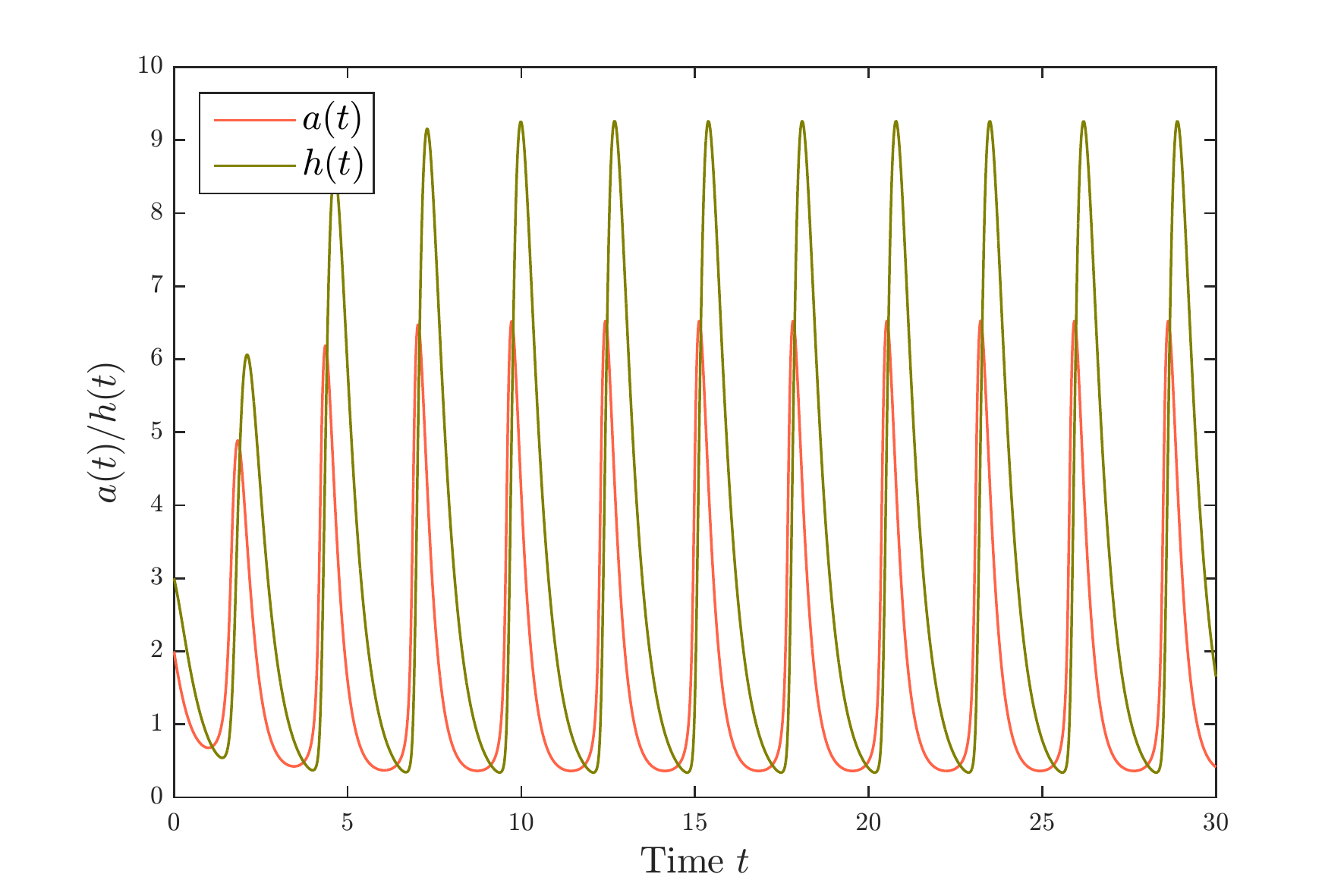}
\caption{The numerical solution of system (\protect \ref{6.9}) for $\protect%
\alpha =1$ with the parameters listed in Table \protect \ref{GM_table_par}.}
\label{GM_model_Fig5_1_Unstable}
\end{figure}

\begin{figure}[tbp]
\centering
\includegraphics[width=11cm]{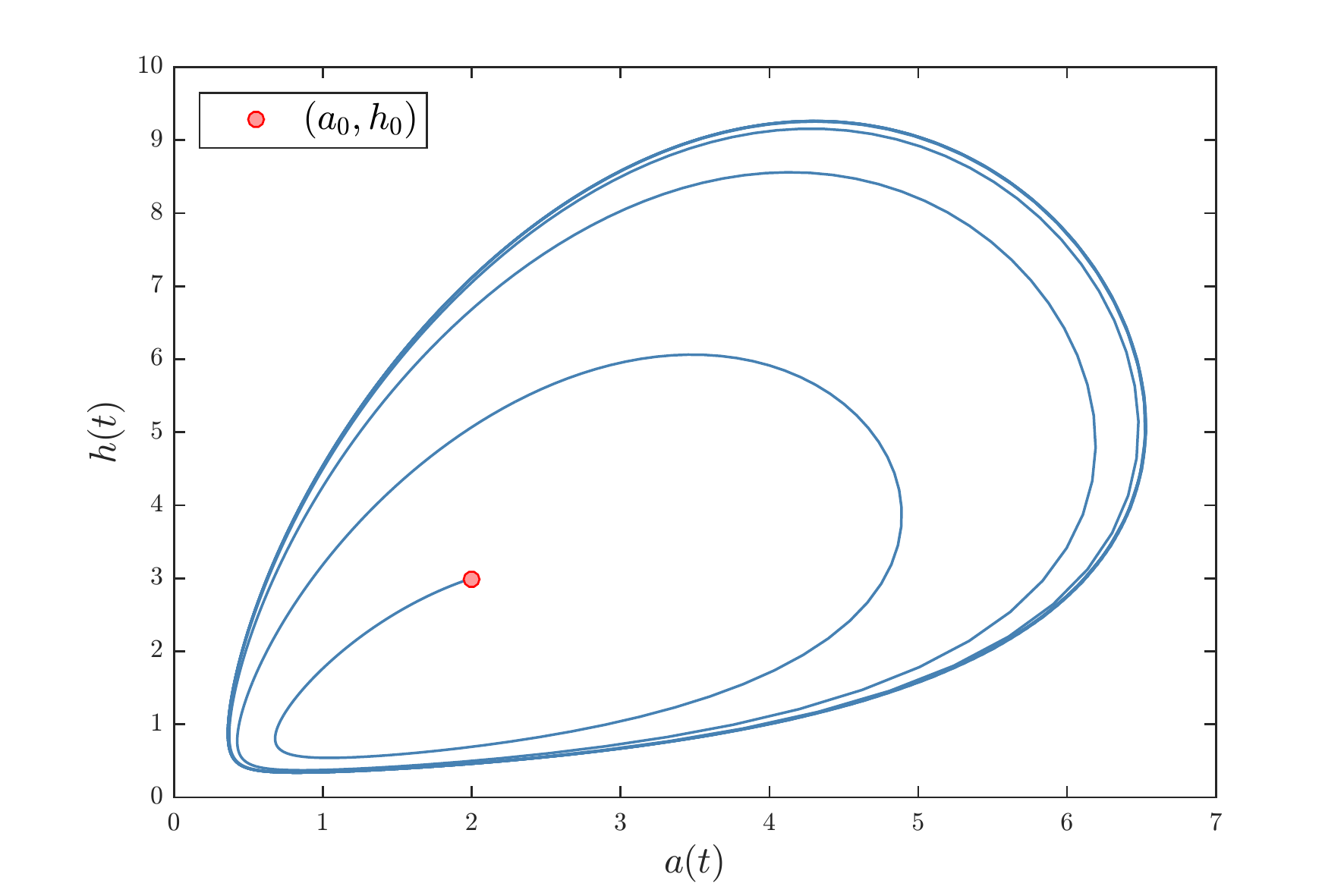}
\caption{Phase plot of system (\protect \ref{6.9}) for $\protect \alpha =1$
with the parameters listed in Table \protect \ref{GM_table_par}.}
\label{GM_model_Fig6_2_phseP_stable}
\end{figure}
\end{example}

\section{Conclusion}
In this paper, we have employed one/two steps first/second order Newton polynomial interpolation to derive new two methods to solve fractional differential equations for several definitions of the fractional derivative, the first one method is the improved version of the Atangana-Seda method which has been widely used in a short time period since its appearance, and the second one method we have proposed new predictor-corrector method and we have used improved Atangana-Seda scheme as a predictor term. The proposed methods have demonstrated their effectiveness with the various examples presented and have proven effective for obtaining accurate approximate solutions for complex systems. The simplicity of displaying proposed methods equations enables us to easily convert them into algorithms and translate them into different programming languages for use in numerical simulations of systems modeling various phenomena in the real world. These methods will open new horizons in the field of numerical analysis of fractional differential equations with many definitions of fractional derivative.

\end{document}